[1994/12/01]

\documentclass[a4paper, 11pt]{amsart}
\usepackage{amssymb,amsxtra}
\usepackage{amsmath,amsthm,amscd}
\usepackage[all,cmtip]{xy}

\chardef\bslash=`\\ 

\newtheorem{thm}{Theorem}[section]
\newtheorem{cor}[thm]{Corollary}
\newtheorem{lem}[thm]{Lemma}
\newtheorem{prop}[thm]{Proposition}

\theoremstyle{definition}
\newtheorem{defn}{Definition}[section]

\theoremstyle{remark}
\newtheorem{rem}{Remark}[section]

\newcommand{\thmref}[1]{Theorem~\ref{#1}}

\newcommand{\lemref}[1]{Lemma~\ref{#1}}
\newcommand{\propref}[1]{Proposition~\ref{#1}}
\newcommand{\corref}[1]{Corollary~\ref{#1}}
\newcommand{\defnref}[1]{Definition~\ref{#1}}

\newcommand{\remref}[1]{Remark~\ref{#1}}

\newcommand{\eval}[2][\right]{\relax
  \ifx#1\right\relax \left.\fi#2#1\rvert}

\def\projlim{\mathop{\oalign{lim\cr
\hidewidth$\longleftarrow$\hidewidth\cr}}}

\def\indlim{\mathop{\oalign{lim\cr
\hidewidth$\longrightarrow$\hidewidth\cr}}}

\def\R{{\mathbf{R}}}
\def\C{{\mathbf{C}}}
\def\Z{{\mathbf{Z}}}
\def\N{{\mathbf{N}}}

\def\t{{\mathbf{T}}}
\def\Q{{\mathbf{Q}}}
\def\sol{{\mathcal{S}}}

\def\defequal{\stackrel{\mathrm{def}}{=}}

\def\<{\langle}
\def\>{\rangle}

\def\tfae{the following conditions are equivalent:}
\def\lan{\langle}
\def\ran{\rangle}

\title{Morita equivalence classification for irrational noncommutative solenoids}

\author[P. Carrillo Rouse]{Paulo Carrillo Rouse}
\address{Institut de Math\'ematiques \\  Université de Toulouse\\
  31062 Toulouse, FRANCE}
\email{paulo.carrillo@math.univ-toulouse.fr}

\author[L. Guillaume]{Laurent Guillaume}
\address{QWANTO\\ 2 rue d'Austerlitz\\ 31000 Toulouse, FRANCE}
\email{laurent@qwanto.fr}

\author{Jean Renault}
\address{Institut Denis Poisson (UMR 7013)\\
  Universit\'e d'Orl\'eans et CNRS \\ 45067
  Orl\'eans Cedex 2, FRANCE}
\email{jean.renault@univ-orleans.fr}

\usepackage{pdfsync}
\begin{document}
\vskip5mm
\begin{abstract}
Given a prime integer $p$, we define a $p$-angle to be an element of the $p$-solenoid ${\mathcal S}_p$ or equivalently an element of $\R\times \Q_p$ modulo $\Z[1/p]$. Noncommutative $p$-solenoids, or equivalently $p$-rotation algebras, are $C^*$-algebras defined by a $p$-angle in the same way as non-commutative tori or usual rotation algebras are defined by an angle. We give a complete classification, in terms of their $p$-angles and in the irrational case, of these $C^*$-algebras up to strong Morita equivalence, but also the associated groupoids up to Morita equivalence; more precisely, we obtain the following classification result:  two irrational $p$-angles define strong Morita equivalent noncommutative solenoids if and only if they define isomorphic dimension groups if and only if they are congruent modulo $GL_2(\Z[1/p])$ if and only if they define Morita equivalent groupoids. This generalises, with a different proof, the classification up to isomorphism, previously obtained by Latr\'emoli\`ere and Packer.
We also give a geometrical construction of projective $C^*$-modules in terms of proper $G$-spaces which exhaust the dimension group up to a shift and implement all relevant Morita equivalences by groupoid equivalences.
\end{abstract}
\maketitle
\markboth{Carrillo-Rouse Guillaume Renault} {NC Solenoids}

\renewcommand{\sectionmark}[1]{}

\section*{Introduction}

The irrational rotation C*-algebras, which have been studied extensively and classified by Rieffel in \cite{rie:irrational} and which appear under various guises such as noncommutative tori or Kronecker foliation C*-algebras, are basic to the theory. We show in this article that, replacing $\Z$ by $\Z[1/p]$, where $p$ is a prime number, the study of the corresponding C*-algebras can be done in exactly the same fashion. Since the dimension group is a complete invariant of these C*-algebras, a large part of the work is devoted to its study. Another feature of the article is the emphasis of geometrical constructions. In particular, Morita equivalence of our C*-algebras are naturally implemented by groupoid equivalences as in \cite{mrw:equivalence} and dimension of projective modules are obtained as volumes. In fact, this geometrical approach was the main motivation for writing this article as a complement to the article \cite{lp:ncsol} by Latr\'emoli\`ere and Packer which introduces these algebras and studies them through C*-algebraic techniques.\\

Here is a description of the content of each section and its main results. The first section gives the prerequisites about the field $\Q_p$ of $p$-adic numbers and its ring of integers $\Z_p$. The main actor is the ring $\Z[1/p]$, whose dual group is the $p$-solenoid $\sol_p$. The key fact is that, just as $\Z$ is a discrete subgroup of $\R$ with quotient $\bf S^1$, $\Z[1/p]$ is a discrete subgroup of $\R\times\Q_p$ with quotient $\sol_p$.\\
The second section is devoted to $p$-angles and their dimension groups. This is just algebra. A $p$-angle is an element of $\sol_p$, or equivalently an element of $\R\times\Q_p$ modulo $\Z[1/p]$. The dimension group of a $p$-angle admits two definitions, either as an inductive limit or, better, as a subgroup of $\Z[1/p]^2$ (\defnref{defn2}). It is interesting there to see the appearance of two kinds of positivity: $\R_+$ in $\R$ and $\Z_p$ in $\Q_p$. The main result, \corref{thm1}, is that the dimension groups are isomorphic if and only if the $p$-angles are congruent modulo ${\rm GL}(2,\Z[1/p])$.\\
The third section gives the definition of the noncommutative solenoids and of the irrational $p$-rotation algebras. In order to show that these C*-algebras are isomorphic and also that they are the universal C*-algebra of canonical commutation relations, we give a general result, \thmref{isomorphisms}, which uses partial Fourier transforms and which can be deduced from \cite[Theorem 7.8]{ren:atge}. It is then easy to exhibit noncommutative solenoids as inductive limits of (usual) irrational algebras, which is a basic result of \cite{lp:ncsol}, and to deduce that their dimension groups are those defined in the previous section.\\
The fourth section gives a geometrical construction of projective C*-modules and of Morita equivalence. It is developed in a rather general setting which includes the usual irrational rotation algebra. Instead of using the C*-algebra of the Kronecker foliation as a stable model for the irrational rotation algebra, we use instead the C*-algebra of the action of $\Z^2$ on $\R$ through the embedding sending $(k,l)$ to $k\alpha+l$ (by \thmref{isomorphisms}, these C*-algebras are isomorphic). Similarly, our stable model is the transformation groupoid $G$ of the action of $\Z[1/p]^2$ on $\R\times\Q_p$ through the embedding defined by the $p$-angle. The first result (\thmref{reduction} and \corref{cor}) gives a groupoid equivalence $Z$ of $G=\Gamma\ltimes X$ and $\Gamma/\Gamma_1\ltimes X/\Gamma_1$ where $\Gamma$ is a locally compact group acting continuously on a locally compact space $X$ and $\Gamma_1$ is a closed normal subgroup of $\Gamma$ which acts freely and properly on $X$. Moreover (\propref{dimension}), if $\mu$ is an invariant measure on $X$, $\Gamma$ is discrete and $X/\Gamma_1$ is compact, then the dimension of $C^*(Z)$, as a projective C*-module over $C^*(G)$ and with respect to the trace defined by $\mu$, is the transverse measure of $X/\Gamma_1$. We apply these results to our action of $\Gamma=\Z[1/p]^2$ on $X=\R\times\Q_p$. By choosing as $\Gamma_1$ a one-dimensional subgroup of $\Gamma$, one exhausts the positive part of the dimension group of $C^*(G)$ of degree 0 (\propref{Psi}). When $\Gamma_1$ is maximal, the groupoid $C^*$-algebra of $\Gamma/\Gamma_1\ltimes X/\Gamma_1$ is a noncommutative $p$-solenoid (\propref{groupoid equivalence}). The last result (\propref{index}) describes the topological structure of the groupoid equivalence between $\Gamma/\Gamma_1\ltimes X/\Gamma_1$ and $\Gamma/\Gamma_2\ltimes X/\Gamma_2$ associated to two subgroups $\Gamma_1$ and $\Gamma_2$ satisfying our assumptions.\\
The fifth section summarises our results about the classification up to Morita equivalence of noncommutative solenoids in a single statement, theorem \ref{thm1grpdsbis}.\\
The appendix reproduces results from \cite{ren:rio} which are used in the text.

\section{Preliminaries}

We recall the definitions and some basic results used in the text. Our reference on the subject is  \cite[Chapter 1]{rob:p-adic}, to which we direct the reader for a detailed exposition. We assume that $p$ is a prime integer although a number of definitions and results are valid when $p$ is an arbitrary integer. The main objects are

$\Z[1/p]=$ group of $p$-adic rational numbers; it is also a Euclidean ring;

${\mathcal S}_p=$ $p$-solenoid;

$\Q_p=$ field of $p$-adic numbers;

$\Z_p=$ ring of $p$-adic integers.

Since $\Z[1/p]$ is a subgroup of $\Q$, it is both a subgroup of $\R$ and of $\Q_p$. It is dense both in $\R$ and in $\Q_p$. Elements of $\Q_p$ are written as Laurent series
\[x=\sum_{n=n_0}^\infty x_np^n,\quad{\rm where}\quad x_n\in\{0,1,\ldots,p-1\}.\]
If $x_{n_0}\not=0$, $n_0=v(x)$ is the valuation of $x$. The $p$-adic norm is $|x|_p=p^{-v(x)}$. The $p$-adic integers are the $p$-adic numbers such that $v(x)\ge 0$. Every $x\in \Q_p$ has an integral part 
\[[x]=\sum_{n=0}^\infty x_np^n\quad\hbox{belonging to}\,\Z_p.\]
One defines its fractional part $\{x\}=x-[x]$ which belongs to $\Z[1/p]$. 

Through the Euclidean algorithm, $\N$ embeds into $\Z_p$. One deduces that $\Q$ embeds into $\Q_p$. It is a dense subfield.
\vskip 3mm
The fundamental character of $\R$ is $e(t)=\exp(i2\pi t)$. It establishes an isomorphism of $\R$ and its dual $\hat\R$: $t\mapsto\chi_t$, where $\chi_t(s)=e(ts)$. Similarly, the fundamental character of $\Q_p$ is $e_p(x)=e(\{x\})$, where $\{x\}$ is the fractional part of $x$. It establishes an isomorphism of $\Q_p$ and its dual $\hat\Q_p$ : $x\mapsto\chi_x$, where $\chi_x(y)=e_p(xy)$.

The $p$-solenoid ${\mathcal S}_p$ is defined as the dual group of $\Z[1/p]$. Since $\Z[1/p]$ is the inductive limit
\[\Z[1/p]=\indlim{(\Z\buildrel{p}\over\longrightarrow\Z\buildrel{p}\over\longrightarrow\ldots)}\]
${\mathcal S}_p$ is the projective limit
\[\sol_p=\projlim(\t \buildrel{z^p}\over\longleftarrow \t \buildrel{z^p}\over\longleftarrow\ldots)\]
where $\t=\{z\in\C: |z|=1\}$ is identified with $\hat\Z$ through the duality $\langle z, n\rangle=z^n$.
\vskip 3mm
The most important results about $\Z[1/p]$ are analogous to classical properties of $\Z$:
\begin{enumerate}
 \item $\Z$ is a discrete subgroup of $\R$;
 \item the quotient group $\R/\Z$ is identified with $\hat\Z$.
\end{enumerate}

\begin{thm}
 $\Z[1/p]$ embeds into $\R\times\Q_p$ as a discrete subgroup.
\end{thm}
\begin{proof}
We choose the embedding $x\mapsto (x,-x)$. For later purposes, we write $j:\Z[1/p]\to\R$ the embedding into $\R$, $j_p:\Z[1/p]\to\Q_p$ the embedding into $Q_p$ and $J=(j,-j_p)$ this embedding into $\R\times\Q_p$. We observe that $V=(-1,1)\times\Z_p$  is an open subset of $\R\times\Q_p$. We show that $J(\Z[1/p])\cap V=\{(0,0)\}$. Indeed, if $x=\sum_m^nx_ip^i\in\Z[1/p]$ belongs to $\Z_p$, then $m\ge 0$. This implies $|x|\ge 1$ unless $x=0$.
\end{proof}
Note that instead of $J=(j,-j_p)$, we could have used $J'=(j,j_p)$.

\begin{thm}
 The quotient group $(\R\times\Q_p)/\Z[1/p]$ is identified with the $p$-solenoid $\sol_p$.
\end{thm}
\begin{proof}
 By duality, $J,J':\Z[1/p]\to\R\times\Q_p$ which are strict and injective gives $\hat J,\hat J':\hat\R\times\hat\Q_p\to\sol_p$ strict and surjective. We identify $\R$ and $\Q_p$ with their duals as above. We check that $\hat j:\R\to\sol_p$ is given by $\hat j(t)=(e(t/p^n))$ and that $\hat j_p:\Q_p\to\sol_p$ is given by $\hat j_p(x)=(e_p(x/p^n))$. Thus $\hat J, \hat J'$ are given by $\hat J(t,x)=(e(t/p^n-\{x/p^n\})$ and $\hat J'(t,x)=(e(t/p^n+\{x/p^n\})$.  They are surjective. One can check that the kernel of $\hat J'$ is the image of $J$.
\end{proof}

We retain the exact sequence
\[0\to\Z[1/p]\buildrel{J}\over\longrightarrow\R\times\Q_p\buildrel{\hat J'}\over\rightarrow{\sol_p}\to 1.\]

\begin{lem}
 Let $H$ be a subgroup of a topological group $G$. Suppose that there exists a neighborhood $V$ of $1\in G$ such that $V\cap H=\{1\}$. Then,
 \begin{enumerate}
 \item $H$ is a discrete subgroup;
 \item every subset $U$ of $G$ which meets every coset $xH$ and such that $U^{-1}U\subset V$ is a fundamental domain for the quotient map $\pi:G\to G/H$.
\end{enumerate}
\end{lem}

\begin{proof}
 For all $h\in H$, $Vh\cap H=\{h\}$. If $U$ meets every coset $xH$, $\pi(U)=G/H$. Suppose that $x,y\in U$ and $\pi(x)=\pi(y)$. Then $x^{-1}y$ belongs to $U^{-1}U\cap H$, hence is equal to 1.
\end{proof}

We apply this lemma to $H=J(\Z[1/p])$ and $G=\R\times\Q_p$. 

\begin{cor} $[0,1)\times \Z_p$ is a fundamental domain for the quotient map 
\[\pi=\hat J':\R\times\Q_p\to\sol_p\] 
sending $(\alpha,{\bf a})$ to $(e(\alpha/p^n+\{{\bf a}/p^n\}))_n$.\end{cor}

\vskip2cm

\section{$p$-angles and their dimension groups}
\vskip3mm
In this section, $\pi:\R\times\Q_p\to\sol_p$ denotes the quotient map $\pi(\alpha,{\bf a})=(e(\alpha/p^n+\{{\bf a}/p^n\}))_n$ and $J:\Z[1/p]\to \R\times\Q_p$ the inclusion $J(x)=(x,-x)$. We have the exact sequence
\[0\to\Z[1/p]\buildrel{J}\over\longrightarrow\R\times\Q_p\buildrel{\pi}\over\rightarrow{\sol_p}\to 1.\]
In analogy with the usual definition of an angle as a real number modulo $\Z$, we give the following definition.

\begin{defn}
 We define a $p$-angle as an element $\zeta=(\zeta_n)$ of $\sol_p$. Equivalently, we view it as an element ${\vec\alpha}=(\alpha,{\bf a})\in\R\times\Q_p$ modulo $\Z[1/p]$: $(\alpha,{\bf a})$ and $(\alpha+x, {\bf a}-x)$ define the same $p$-angle for all $x\in\Z[1/p]$. The identification is given by $\zeta=\pi(\alpha,{\bf a})$, i.e. $\zeta_n=e(\alpha/p^n+\{{\bf a}/p^n\})$ for all $n\in\N$.
  \end{defn}

Note that there is one and only one representative $(\alpha,{\bf a})$ of a $p$-angle $\zeta$ in the fundamental domain $[0,1)\times\Z_p$.  If $(\alpha,{\bf a})\in [0,1)\times\Z_p$, we write 
 \[{\bf a}=\sum_{i=0}^\infty a_ip^i\quad{\rm where}\quad a_i\in\{0,1,\ldots, p-1\}\]
 and decompose ${\bf a}/p^n$ as the sum of its fractional and its integer parts:
 \[\dfrac{{\bf a}}{p^n}=\{\dfrac{{\bf a}}{p^n}\}+[\dfrac{{\bf a}}{p^n}].\]
 
We then introduce the following notations:
 
\begin{eqnarray*}
 \alpha_n=&\dfrac{\alpha}{p^n}+\{\dfrac{{\bf a}}{p^n}\}=&\dfrac{1}{p^n}(\alpha+\sum_{i=0}^{n-1}a_ip^i)\cr
 {\bf a}_n=&\dfrac{{\bf a}}{p^n}-\{\dfrac{{\bf a}}{p^n}\}=&\dfrac{1}{p^n}\sum_{i=n}^\infty a_ip^i=[\dfrac{{\bf a}}{p^n}]\cr
 {\bf a}^{(n)}=&p^n\{\dfrac{{\bf a}}{p^n}\}=&\sum_{i=0}^{n-1}a_ip^i
\end{eqnarray*}
We note the induction relation
\[p\alpha_{n+1}=\alpha_n+a_n.\]

\begin{lem}
 Using above notations, we have: 
\begin{enumerate}
\item $(\alpha_n,{\bf a}_n)$ define the same $p$-angle as $p^{-n}(\alpha,{\bf a})$,  
\item $\zeta_n=e(\alpha_n)$,
 \item if $0<\alpha<1$, then $0<\alpha_n<1$ for all $n\in\N$,
 \item $p\alpha_{n+1}=\alpha_n+a_n$.
 \end{enumerate}
\end{lem}

Note that it is possible, as Latr\'emoli\`ere and Packer do, to define a $p$-angle as a sequence $(\alpha_n)_{n\in\N}$ of real numbers such that $0\le\alpha_0<1$ and $p\alpha_{n+1}-\alpha_n$ belongs to $\{0,1,\ldots,p-1\}$. It defines $\zeta\in\sol_p$ such that $\zeta_n=e(\alpha_n)$. One recovers $a_n=p\alpha_{n+1}-\alpha_n$ and ${\bf a}=\sum_{i=0}^\infty a_ip^i$ in $\Z_p$.

\begin{lem}
Let $\zeta=(\zeta_n)\in \sol_p$ be a $p$-angle given by $(\alpha,{\bf a})\in (0,1)\times\Z_p$, then
 
\begin{enumerate}
 \item  $\overline\zeta=(\overline\zeta_n)$ is a $p$-angle which is given by $(\beta,{\bf b})\in (0,1)\times\Z_p$ where $\beta=1-\alpha$ and ${\bf b}=\sigma({\bf a})\defequal\sum_{i=0}^\infty(p-1-a_i)p^i$. Moreover, we have then $\beta_n=1-\alpha_n$ for all $n\in\N$.
 \item the $k$-shift $\eta$ of $\zeta$, denoted by $\eta=S_k\zeta$ and defined by $\eta_n=\zeta_{n+k}$, is a $p$-angle. It is given by $(\alpha_k,{\bf a}_k)$, or equivalently by $p^{-k}(\alpha,{\bf a})$.
\end{enumerate}

\end{lem}

\begin{proof}
 For (i), one uses the formula $-{\bf a}=1+\sigma({\bf a})$. The assertion (ii) is clear.
\end{proof}

\begin{defn}
 We then say that
 
\begin{enumerate}
 \item $(1-\alpha,\sigma({\bf a}))$ is the conjugate of $(\alpha,{\bf a})$,
 \item $(\alpha_k,{\bf a}_k)$ is the $k$-shift of $(\alpha,{\bf a})$. We shall sometimes write: $\alpha_k=S_k\alpha$, ${\bf a}_k=S_k{\bf a}$ and $S_k(\alpha,{\bf a})=(\alpha_k, {\bf a}_k)$.
\end{enumerate}
\end{defn}

As seen earlier, in terms of $(\alpha,{\bf a})$,  $k$-shifting amounts to multiplying by $p^{-k}$. Indeed,

\[\vec\alpha_k\defequal(\alpha_k,{\bf a}_k)=p^{-k}(\alpha+\sum_{i=0}^{n-1}a_ip^i,{\bf a}-\sum_{i=0}^{n-1}a_ip^i)\]
is equivalent to $p^{-k}\vec\alpha$. In terms of $\zeta$, $k$-shifting amounts to raising at the power $p^{-k}$.
We may write $S_k\zeta=\zeta^{(p^{-k})}$.\\

Both operations, conjugation and shifts, can be expressed under a matricial form. Given a matrix $\left(\begin{matrix}
 a&c\\
 b&d
\end{matrix}\right)$ with coefficients in $\Z[1/p]$ and $x$ in $\R$ or $\Q_p$ such that $bx+d\not=0$ , we write:
\[\left(\begin{matrix}
 a&c\\
 b&d
\end{matrix}\right)x=\dfrac{ax+c}{bx+d}.\]
Given $(s,x)\in\R\times\Q_p$, we write 
\[(t,y)=\left(\begin{matrix}
 a&c\\
 b&d
\end{matrix}\right)(s,x)\quad{\rm if}\quad t=\left(\begin{matrix}
 a&c\\
 b&d
\end{matrix}\right)s\quad{\rm and}\quad y=\left(\begin{matrix}
 a&-c\\
 -b&d
\end{matrix}\right)x.\]

\begin{defn}\label{defcongruent}
 We say that two elements $(\alpha, {\bf a})$ and $(\beta, {\bf b})$ of $\R\times\Q_p$ are
 
\begin{enumerate}
\item congruent through the matrix $\left(\begin{matrix}
 a&c\\
 b&d
\end{matrix}\right)\in{\rm GL}(2,\Z[1/p])$ if 
\[(\alpha, {\bf a})=\left(\begin{matrix}
 a&c\\
 b&d
\end{matrix}\right)(\beta, {\bf b}),\]
\item congruent modulo ${\rm GL}(2,\Z[1/p])$[resp. modulo ${\rm GL}(2,\Z)$] if they are congruent through a matrix in ${\rm GL}(2,\Z[1/p])$ [resp. in ${\rm GL}(2,\Z)$].
\end{enumerate}
\end{defn}

\begin{rem}
 Since the centre 
 \[Z=\{\left(\begin{matrix}
 \epsilon p^{-n}&0\\
 0& \epsilon p^{-n}
\end{matrix}\right),n\in\Z, \epsilon=\pm 1\}\]
 of ${\rm GL}(2,\Z[1/p])$ acts as the identity, this action of  ${\rm GL}(2,\Z[1/p])$ on $\R\times\Q_p$ factors through ${\rm PGL}(2,\Z[1/p])={\rm GL}(2,\Z[1/p])/Z$. Thus we could say congruent modulo ${\rm PGL}(2,\Z[1/p])$ instead of ${\rm GL}(2,\Z[1/p])$.
\end{rem}

\begin{lem} \label{positivity} If $(\alpha, {\bf a})$ and $(\beta, {\bf b})$ are congruent modulo ${\rm GL}(2,\Z[1/p])$, there exists $\left(\begin{matrix}
 a&c\\
 b&d
\end{matrix}\right)\in {\rm GL}(2,\Z[1/p])$ such that
\[(\alpha, {\bf a})=\left(\begin{matrix}
 a&c\\
 b&d
\end{matrix}\right)(\beta, {\bf b})\]
and satisfying $b\beta +d>0$ and $v(-b{\bf b}+d)=0$.
\end{lem}

\begin{proof}
 By assumption, there exists $\left(\begin{matrix}
 a'&c'\\
 b'&d'
\end{matrix}\right)\in {\rm GL}(2,\Z[1/p])$ such that
\[(\alpha, {\bf a})=\left(\begin{matrix}
 a'&c'\\
 b'&d'
\end{matrix}\right)(\beta, {\bf b}).\]
We let $\epsilon$ be the sign of $b'\beta +d'$ and $n=v(-b'{\bf b}+d')$. Then
\[\left(\begin{matrix}
 a&c\\
 b&d
\end{matrix}\right)=\left(\begin{matrix}
 \epsilon p^{-n}&0\\
 0& \epsilon p^{-n}
\end{matrix}\right)\left(\begin{matrix}
 a'&c'\\
 b'&d'
\end{matrix}\right)\]
satisfies the required conditions.
\end{proof}

\begin{prop}\label{propcongruent}
Let $(\alpha, {\bf a})\in \R\times\Q_p$.
 \begin{enumerate}
 \item given $x\in\Z[1/p]$, $(\alpha,{\bf a})$ and $(\alpha+x,{\bf a}-x)$ are congruent through the matrix $\left(\begin{matrix}
 1&x\\
 0&1
\end{matrix}\right)$.
\item $(1-\alpha,\sigma({\bf a}))$ and $(\alpha,{\bf a})$ are congruent through the matrix  $\left(\begin{matrix}
 -1&1\\
 0&1
\end{matrix}\right)$.
\item given $k\in\N$,  $S_k(\alpha,{\bf a})$ and $(\alpha,{\bf a})$ are congruent through the matrix $\left(\begin{matrix}
  p^{-k}&\{{p^{-k}\bf a}\}\\
 0&1
\end{matrix}\right)$.
\end{enumerate}
\end{prop}
\vskip5mm

Before defining the dimension group of a $p$-angle, let us review the case of the irrational rotation.
By definition, the dimension group of an irrational real number $\alpha$ is the subgroup $\Z\alpha+\Z$ of $\R$ ordered by the usual order of $\R$. The classification of these dimension groups is recalled in the next theorem, along with its algebraic proof, because we shall have a similar result with a similar proof for the dimension groups of $p$-angles. We note however that Rieffel gives in \cite{rie:irrational} an interesting proof using the Morita equivalence of the corresponding C*-algebras instead of \thmref{r}.

\begin{thm}\cite[Theorem 14.23]{goo:poag}\label{rot}  Let $A$ denote the set of positive irrational real numbers.
Let $\alpha$ and $\beta$ be non-zero irrational real numbers. Then $\Z\alpha+\Z$ and $\Z\beta+\Z$ are isomorphic as ordered groups if and only if there exist  $a,b,c,d$ in $\Z$ such that
 \begin{enumerate}
 \item  $ad-bc=\pm 1$,
 \item  $b\beta+d>0$
 \item $\alpha=(a\beta+c)/(b\beta+d)$.
 \end{enumerate}
Moreover, $(\Z\alpha+\Z,1)$ and $(\Z\beta+\Z,1)$ are isomorphic as ordered groups with units if and only if $b=0$. In this case, $d=1$ and $a=\pm 1$.
\end{thm}

\begin{proof} The proof is deduced from the following diagram:
 
 \[\begin{CD}
\Z\times\Z@>{j_\alpha}>{}>\Z\alpha+\Z\\
 @V{M}V{}V  @V{\times r}V{}V\\
\Z\times\Z @>{j_\beta}>{}>\Z\beta+\Z
\end{CD}\]
where $M=\left(\begin{matrix}
 a&b\\
 c&d
\end{matrix}\right)$ and $j_\alpha(x,y)=x\alpha+y$. If $M$ exists, one obtains a commutative diagram by setting $r=b\beta+d$, hence an order isomorphism from $\Z\alpha+\Z$ onto $\Z\beta+\Z$. Conversely, suppose that $\varphi:\Z\alpha+\Z\to\Z\beta+\Z$ is an isomorphism of ordered subgroups. As a group isomorphism, it is given by a matrix $M=\left(\begin{matrix}
 a&b\\
 c&d
\end{matrix}\right)$ in ${\rm GL}(2,\Z)$. Moreover, according to a result by Hion recalled in \thmref{r}, an order isomorphism of two subgroups of $\R$ is necessarily the multiplication by some $r>0$. This gives the above commuting diagram. Necessarily $r\alpha=a\beta+c$ and  $r=r1=b\beta+d$. Alternatively, $a,b,c,d$ can be obtained from these relations. In the unit preserving case, $r=1$. Therefore $b=0$ and $d=1$.
\end{proof}

We give now the definition of the dimension group of a $p$-angle as an inductive limit. An intrinsic equivalent definition will be given afterwards (\defnref{defn2}).

\begin{defn}\label{defn1} Let $(\alpha,{\bf a})\in\R\times\Z_p$. We define the dimension group $G(\alpha,{\bf a})$ as the subgroup of $\R$ generated by $1$ and the sequence $(\alpha_n)_{n\in\N}$, ordered by the usual order of $\R$.
\end{defn}

\begin{prop}Let $(\alpha,{\bf a})\in\R\times\Z_p$ and $G_n=\Z\alpha_n+\Z$ for $n\in\N$.
\begin{enumerate}
 \item the sequence of subgroups $(G_n)_{n\in\N}$ is increasing and $G(\alpha,{\bf a})=\bigcup_{n\in\N} G_n$.
 \item for all $n\in\N$, $G(\alpha_n,{\bf a}_n)=G(\alpha,{\bf a})$.
 \end{enumerate}
 \end{prop} 
\begin{proof}
 (i) Since $\alpha_n=p\alpha_{n+1}-a_n$, it belongs to $\Z\alpha_{n+1}+\Z$. Therefore, $\Z\alpha_n+\Z\subset \Z\alpha_{n+1}+\Z$.\\
 (ii) $G(\alpha_n,{\bf a}_n)=\bigcup_{i\ge n} G_i=\bigcup_{i\ge 0} G_i=G(\alpha,{\bf a})$.
\end{proof}
\begin{cor}
 The dimension group $G(\alpha,{\bf a})$ depends only on the $p$-angle $[\alpha,{\bf a}]$ defined by $(\alpha,{\bf a})$. From now on, it will be denoted by $G[\alpha,{\bf a}]$.
\end{cor}
Our study will be limited to irrational $p$-angles as defined below.

\begin{lem} For $\zeta=\pi(\alpha,{\bf a})\in{\mathcal S}_p$, \tfae
\begin{enumerate}
 \item there exists $n\in\N$ such that $\zeta_n$ is a root of unity;
 \item there exists $n\in\N$ such that $\alpha_n$ is rational;
 \item for all $n\in\N$, $\zeta_n$ is a root of unity;
 \item for all $n\in\N$, $\alpha_n$ is rational.
\end{enumerate}
 \end{lem}
 
 \begin{proof}
 Suppose that $\zeta_n^k=1$ where $k\in\N\setminus\{0\}$. Then $\zeta_{n-1}^k=\zeta_n^{pk}=1$. On the other hand, the relation $\zeta_n=\zeta_{n+1}^p$ gives $\zeta_{n+1}^{pk}=\zeta_n^k=1$.
\end{proof}

\vskip2mm 

\begin{defn}
 We say that $\zeta\in {\mathcal S}_p$ is a root of unity or that $(\alpha_0,{\bf a}_0)\in \R\times\Q_p$ is rational  if the conditions of the lemma are satisfied and irrational otherwise. \end{defn}
 
 If $(\alpha,{\bf a})$ is irrational, we identify $G_n=\Z\alpha_n+\Z$ and $\Z^2$ through the basis $(\alpha_n,1)$. The matrix of the embedding of $G_n$ into $G_{n+1}$ is $\left(\begin{matrix}
p&0\\
-a_n&1
\end{matrix}\right)$. The matrix of the embedding of $G_0$ into $G_n$ is $\left(\begin{matrix}
p^n&0\\
-{\bf a}^{(n)}&1
\end{matrix}\right)$.

\vskip 5mm 

Our second definition of the dimension group of an element of $\R\times\Q_p$ is based on the following observation.

\begin{prop} \label{second}Let $(\alpha,{\bf a})\in\R\times\Z_p$. The map from $G[\alpha,{\bf a}]$ to $\Z[1/p]\times\Z[1/p]$ sending 
 \[k\alpha_n+l\quad\hbox{to}\quad(p^{-n}k, l+k\{p^{-n}{\bf a}\} )\]
 \begin{enumerate}
 \item is well defined;
 \item is a group homomorphism;
 \item injective;
 \item its range consists of the elements $(x,y)\in \Z[1/p]\times\Z[1/p]$ such that $x{\bf a}-y\in\Z_p$;
 \item its inverse sends $(x,y)$ to $x\alpha+y$.
 \end{enumerate}
\end{prop}
\begin{proof}
 (i) The map is well defined: if we replace $\alpha_n$ by $p\alpha_{n+1}-a_n$, we have:
 \[k\alpha_n+l=k(p\alpha_{n+1}-a_n)+l=kp\alpha_{n+1}+l-ka_n\]
 which is sent to 
 \[(p^{-(n+1)}kp, l-ka_n+kp\{p^{-(n+1)}{\bf a}\})=(p^{-n}k, l+k\{p^{-n}{\bf a}\} ).\]
 (ii) By choosing $n$ large enough, any pair of elements of $G(\alpha,{\bf a})$ is of the form $(k\alpha_n+l,k'\alpha_n+l')$. Then its sum $(k+k')\alpha_n+(l+l')$ is sent to $(p^{-n}(k+k'), (l+l')+(k+k')\{p^{-n}{\bf a}\})$, which is the sum of the images of $k\alpha_n+l$ and $k'\alpha_n+l'$.\\
 (iii) If $(p^{-n}k, l+k\{p^{-n}{\bf a}\} )=(0,0)$, then $k=0$ and $l=0$. Therefore $k\alpha_n+l=0$.\\
 Suppose that $x=p^{-n}k$ and $y=l+k\{p^{-n}{\bf a}\}$ where $n\in\N$ and $k,l\in\Z$. Then,
 \[x{\bf a}-y=p^{-n}k{\bf a}-l+k\{p^{-n}{\bf a}=k[p^{-n}{\bf a}]-l\]
 belongs to $\Z_p$. Conversely, let $(x,y)\in\Z[1/p]\times\Z[1/p]$ such that $v(x{\bf a}-y)\ge 0$. We have to find $n\in\N$ and $k,l\in\Z$ such that $x=p^{-n}k$ and $y=l+k\{p^{-n}{\bf a}\}$. If $x\in\Z$, we can take $n=0$, $k=x$ and $l=y$. Indeed, since $y-x{\bf a}$ and $x{\bf a}$ belong to $\Z_p$, so does $y=y-x{\bf a}+x{\bf a}$. Since $y$ also belongs to $\Z[1/p]$, $y$ belongs to $\Z$. If $x\in\Z[1/p]\setminus\Z$, we write $x=p^{-n}k$ with $k$ prime to $p$. This determines $n\in\N$ and $k\in\Z$ uniquely. It remains to show that $l=y-k\{p^{-n}{\bf a}\}$ is an integer. Since $l$ belongs to $\Z[1/p]$, it suffices to show that it belongs to $\Z_p$. We have:
\begin{eqnarray*}
l&=&y-k\{p^{-n}{\bf a}\} \cr
&=&y-x{\bf a}+x{\bf a}-k\{p^{-n}{\bf a}\} \cr
&=&y-x{\bf a}+p^{-n}k{\bf a}-k\{p^{-n}{\bf a}\} \cr
&=&y-x{\bf a}+k(p^{-n}{\bf a}-\{p^{-n}{\bf a}\}) \cr
&=&(y-x{\bf a})+k[p^{-n}{\bf a}] \cr
\end{eqnarray*}
Since $y-x{\bf a}$, $k$ and $[p^{-n}{\bf a}]$ belong to $\Z_p$, so does $l$.\\
(iv) If $(x,y)=(p^{-n}k,l+k\{p^{-n}{\bf a})\}$, then
\[x\alpha+y= p^{-n}k\alpha+k\{p^{-n}{\bf a}\}+l=k\alpha_n+l.\]
\end{proof}
\vskip5mm
This proposition motivates our second definition of the dimension group of an element of $\R\times\Q_p$.
\begin{defn}\label{defn2} We define the dimension group of $(\alpha,{\bf a})\in \R\times\Q_p$ as the ordered abelian group
\[G'(\alpha,{\bf a})=\{(x,y)\in\Z[1/p]\times\Z[1/p]: v(x{\bf a}-y)\ge 0\},\]
\[G'(\alpha,{\bf a})^+=\{(x,y)\in G'(\alpha,{\bf a}): x\alpha+y\ge 0\}.\]
When necessary, we shall normalise it by the unit $(0,1)$. 
 \end{defn}
 
 Note that, as a group, $G'(\alpha,{\bf a})$ depends only on $\bf a$.

\begin{prop}\label{isomorphism} Assume that $(\alpha,{\bf a})\in \R\times\Z_p$. Then the map $j_\alpha:(x,y)\mapsto x\alpha+y$ is a unit preserving ordered group isomorphism from $(G'(\alpha,{\bf a}),(0,1))$ onto $(G[\alpha,{\bf a}],1)$.
 \end{prop}

\begin{proof} This is immediate from the definition of $G'(\alpha,{\bf a})$ and \propref{second}.
 
\end{proof}

We recall that the map $x\mapsto\{x\}+\Z$ is a group homomorphism of $\Q_p$ onto $\Z[1/p]/\Z$ having $\Z_p$ as kernel. This gives a convenient description of the quotient $(\Z[1/p]\times\Z[1/p])/G'(\alpha,{\bf a})$:

\begin{prop}\label{tau} Let $(\alpha,{\bf a})\in \R\times\Q_p$. The short sequence of groups
\[0\to G'(\alpha,{\bf a})\buildrel{i}\over\longrightarrow\Z[1/p]\times\Z[1/p]\buildrel{\tau_{\bf a}}\over\rightarrow\Z[1/p]/\Z\to 0,\]
where $i$ is the inclusion and $\tau_{\bf a}(x,y)=\{x{\bf a}-y\}+\Z$, is exact.
\end{prop}

\begin{proof}
 The injectivity of $i$ and the surjectivity of $\tau_{\bf a}$ are clear. By definition $i(G'(\alpha,{\bf a}))$ is the kernel of $\tau_{\bf a}$.
\end{proof}

Contrary to $G[\alpha,{\bf a}]$, $G'(\alpha,{\bf a})$ depends on the representative $(\alpha,{\bf a})$ of a $p$-angle.

\begin{prop}\label{translation} Let $c\in\Z[1/p]$. Then the automorphism $\left(\begin{matrix}
 1&0\\
 c&1
\end{matrix}\right)$ carries $G'(\alpha,{\bf a})$ onto $G'(\alpha-c,{\bf a}+c)$.
\end{prop}

\begin{proof}
 If $x'=x$ and $y'=cx+y$, then $x'(\alpha-c)+y'=x\alpha+y$ and $x'({\bf a}-c)-y'=x{\bf a}-y$.
\end{proof}

The following classical and curiously similar results are the key ingredients of our main theorem.

\begin{thm}(Hion)\label{r} Let $G$ and $H$ be subgroups of $\R$ endowed with the natural order. Assume that $G\not=0$. If $\varphi$ is an order isomorphism from $G$ into $H$, there exists a real number $r> 0$ such that $\varphi(g)=rg$ for all $g\in G$.
\end{thm}

\begin{proof}
 See \cite[Chapter IV, Section 1]{fuc:po} for an elegant proof.
\end{proof}

\begin{thm}\label{rho} Let $f$ be an automorphism of the group $\Z[1/p]/\Z$. There exists  a unit element $\rho$  in $\Z_p$ such that $f(x+\Z)=\{\rho x\}+\Z$ for all $x\in\Z[1/p]$. We write $f=\rho\cdot$. The correspondence $\rho\mapsto\rho\cdot$ is an isomorphism of the group of units $\Z_p^\times$ and ${\rm Aut}(\Z[1/p]/\Z)$. 
\end{thm}

\begin{proof}
 See \cite[Chapter VIII, Section 43, Example 3]{fuc:iagI}.
\end{proof}

\begin{lem}\label{extension}
 Let $G$ be a group and $\varphi:G'(\alpha,{\bf a})\to G$ be a group homomorphism. Then $\varphi$ extends uniquely to a group homomorphism $\tilde\varphi:\Z[1/p]\times\Z[1/p]\to G$.
\end{lem}

\begin{proof}
 We first observe that $\Z[1/p]\times\Z[1/p]=\bigcup_{n\in\N}p^{-n}G'(\alpha,{\bf a})$. Indeed, let $(x,y)\in \Z[1/p]\times\Z[1/p]$ and $m=v(x{\bf a}-y)$. If $m\ge 0$, $(x,y)$ belongs to $G'(\alpha,{\bf a})$ by definition. If $m=-n<0$, $v(p^n(x{\bf a}-y))=n+m=0$, Hence $(x',y')=p^n(x,y)$ belongs to $G'(\alpha,{\bf a})$ and $(x,y)$ belongs to $p^{-n}G'(\alpha,{\bf a})$. Let us show that, if $(x,y)=p^{-n}(x',y')$ with $(x',y')\in G'(\alpha,{\bf a})$, $\tilde\varphi(x,y)=p^{-n}\varphi(x',y')$ defines $\tilde\varphi$ unambiguously. Indeed, if $p^{-n}(x',y')=p^{-m}(x'',y'')$ and $m\le n$, $(x',y')=p^{n-m}(x'',y'')$ and $\varphi(x',y')=p^{n-m}\varphi(x'',y'')$. Therefore $p^{-n}\varphi(x',y')=p^{-m}\varphi(x'',y'')$. Since the restriction of $\tilde\varphi$ to $p^{-n}G'(\alpha,{\bf a})$ is a group homomorphism  and $\Z[1/p]\times\Z[1/p]$ is the increasing union of the $p^{-n}G'(\alpha,{\bf a})$'s, $\tilde\varphi$ is a group homomorphism. The uniqueness is clear.
 \end{proof}

\begin{thm}\label{main} Let $(\alpha,{\bf a})$ and $(\beta,{\bf b})$ in $(\R\setminus\Q)\times \Q_p$.\\
1. Every order group isomorphism $\varphi: G'(\alpha,{\bf a})\to G'(\beta,{\bf b})$ extends uniquely to a group isomorphism $\tilde\varphi:\Z[1/p]\times\Z[1/p]\to\Z[1/p]\times\Z[1/p]$ given by a matrix $\left(\begin{matrix}
 a&b\\
 c&d
\end{matrix}\right)\in {\rm GL}(2,\Z[1/p])$.\\
2. A group automorphism $\psi:\Z[1/p]\times\Z[1/p]\to\Z[1/p]\times\Z[1/p]$ defined by the matrix $\left(\begin{matrix}
 a&b\\
 c&d
\end{matrix}\right)$ restricts to an ordered group isomorphism $\varphi: G'(\alpha,{\bf a})\to G'(\beta,{\bf b})$ if and only if the following four conditions are satisfied
\begin{enumerate}
 \item $v(-b{\bf b}+d)=0$
 \item $(-b{\bf b}+d){\bf a}=a{\bf b}-c$
 \item $b\beta+d>0$
 \item $(b\beta+d)\alpha=a\beta+c$.
\end{enumerate}
 \end{thm}

\begin{proof}  1.  \lemref{extension} gives the existence and the uniqueness of the group isomorphism $\tilde\varphi$. It is known that every group automorphism of $\Z[1/p]\times\Z[1/p]$ is given by a matrix $\left(\begin{matrix}
 a&b\\
 c&d
\end{matrix}\right)\in {\rm GL}(2,\Z[1/p])$.\\
2. Let $\psi$ be an automorphism of $\Z[1/p]\times\Z[1/p]$ given by a matrix $\left(\begin{matrix}
 a&b\\
 c&d
\end{matrix}\right)\in {\rm GL}(2,\Z[1/p])$. Suppose that $\psi(G'(\alpha,{\bf a})=G'(\beta,{\bf b})$. Then $\psi$ descends to an automorphism $\dot\psi$ of $\Z[1/p]/\Z$. According to \thmref{rho}, there exists $\rho\in\Z_p^\times$ such that $\dot\psi=\rho\cdot$. The following diagram, where $\varphi$ is the restriction of $\psi$ to $G'(\alpha,{\bf a})$ and $\tau_{\bf a}, \tau_{\bf b}$ have been defined in \propref{tau}, is commutative.
\[\begin{CD}
G(\alpha,{\bf a})@>{}>{}>\Z[1/p]\times\Z[1/p]@>{\tau_{\bf a}}>{}>\Z[1/p]/\Z\\
 @V{\varphi}V{}V @V{\psi}V{}V  @V{\rho\cdot}V{}V\\
G(\beta,{\bf b})@>{}>{}>
\Z[1/p]\times\Z[1/p] @>{\tau_{\bf b}}>{}>\Z[1/p]/\Z
\end{CD}\]
The commutativity of the diagram gives the equation
\[\{x(\rho{\bf a})-y\}=\{x(a{\bf b}-c)-y(-b{\bf b}+d)\}\quad ({\rm mod}\, 1)\]
for all $x,y\in\Z[1/p]$. By continuity of the fundamental character $e_p(x)=e(\{x\})$ of $\Q_p$ and density of $\Z[1/p]$ in $Q_p$, one obtains $e_p(\rho y)=e_p((-b{\bf b}+d)y)$ for all $y\in\Q_p$. Since $e_p$ establishes a self-duality of $\Q_p$, $\rho=-b{\bf b}+d$. Similarly, one obtains $\rho{\bf a}=a{\bf b}-c$. Suppose furthermore that $\psi$ preserves the orders of $G'(\alpha,{\bf a})$ and $G'(\beta,{\bf b})$. We let $J_\alpha(x,y)=x\alpha+y$ and $J_\beta(x,y)=x\beta+y$ for $(x,y)\in\Z[1/p]\times\Z[1/p]$. Then $\psi$ induces an order isomorphism of $J_\alpha(G'(\alpha,{\bf a}))$ and $J_\beta(G'(\beta,{\bf b}))$ which are subgroups of $\R$. By \thmref{r}, there exists a real number $r>0$ making commutative the following diagram
 \[\begin{CD}
G(\alpha,{\bf a})@>{}>{}>\Z[1/p]\times\Z[1/p]@>{J_\alpha}>{}>\R\\
 @V{\varphi}V{}V @V{\psi}V{}V  @V{r}V{}V\\
G(\beta,{\bf b})@>{}>{}>
\Z[1/p]\times\Z[1/p] @>{J_\beta}>{}>\R
\end{CD}.\]
The commutativity of the diagram gives $r=b\beta+d$ and $r\alpha=a\beta+c$.\\

Conversely, suppose that the conditions (i) to (iv) are satisfied. Given $(x,y)\in\Z[1/p]\times\Z[1/p]$ and $(x',y')=\psi(x,y)$, we have
\[x'{\bf b}-y'=x(a{\bf b}-c)+y(b{\bf b}-d)=(-b{\bf b}+d)(x{\bf a}-y).\]
Therefore, $v(x'{\bf b}-y')=v(x{\bf a}-y)$ if $v(-b{\bf b}+d)=0$. This implies that $\psi(G'(\alpha,{\bf a}))=G'(\beta,{\bf b})$. We also have
\[x'\beta+y'=x(a\beta+c)+y(b\beta+d)=(b\beta+d)(x\alpha+y).\]
Therefore $x'\beta+y'$ has the same sign as $x\alpha+y$ if $b\beta+d>0$. This implies that $\psi(G'(\alpha,{\bf a})^+)=G'(\beta,{\bf b})^+$
 
\end{proof}

\begin{cor}\label{thm1} Let $(\alpha,{\bf a})$ and $(\beta,{\bf b})$ in $(\R\setminus\Q)\times\Q_p$. Then \tfae
\begin{enumerate}
 \item $(\alpha,{\bf a})$ and $(\beta,{\bf b})$ are congruent modulo ${\rm GL}(2,\Z[1/p])$ and
 \item $G'(\alpha,{\bf a})$ and $G'(\beta,{\bf b})$ are isomorphic as ordered groups. 
\end{enumerate}
\end{cor}

\begin{proof} (i)$\Rightarrow$ (ii) .
 According to \lemref{positivity}, there exists $\left(\begin{matrix}
 a&c\\
 b&d
\end{matrix}\right)\in {\rm GL}(2,\Z[1/p])$ such that
\[(\alpha, {\bf a})=\left(\begin{matrix}
 a&c\\
 b&d
\end{matrix}\right)(\beta, {\bf b})\]
and satisfying $b\beta +d>0$ and $v(-b{\bf b}+d)=0$. Then $\psi:\Z[1/p]\times\Z[1/p] \to \Z[1/p]\times\Z[1/p]$ defined by the matrix $\left(\begin{matrix}
 a&b\\
 c&d
\end{matrix}\right)$ satisfies the conditions (i) to (iv) of \thmref{main}. Therefore, it restricts to an ordered group isomorphism $\varphi:G'(\alpha,{\bf a})\to G'(\beta,{\bf b})$.\\

(ii) $\Rightarrow$ (i) Let $\varphi:G'(\alpha,{\bf a})\to G'(\beta,{\bf b})$ is an ordered group isomorphism. Conditions (ii) and (iv) of \thmref{main} says that $(\alpha,{\bf a})$ and $(\beta,{\bf b})$ are congruent through the transpose of the matrix of the extension $\tilde\varphi$.
\end{proof}

\begin{rem}\label{shift} One can show, for example by using the definition of the dimension group as an inductive limit, that the condition (i) can be replaced by the existence of $k,l\in\N$ such that $S_k(\alpha, {\bf a})$ and $S_l(\beta, {\bf b})$ are congruent modulo ${\rm GL}(2,\Z)$.
\end{rem}

\begin{rem}
 This corollary gives the classification of the irrational noncommutative solenoids up to Morita equivalence. 
 \end{rem}

\begin{cor} Let $(\alpha,{\bf a})$ and $(\beta,{\bf b})$ in $(\R\setminus\Q)\times\Q_p$. Then \tfae
\begin{enumerate}
 \item  there exist $k,l\in\N$ such that $S_k(\alpha, {\bf a})$ and $S_l(\beta, {\bf b})$ are either equal or conjugate;
  \item the dimension groups $(G'(\alpha,{\bf a}),u)$ and $(G'(\beta,{\bf b}),u)$ are isomorphic as ordered groups with units, where the unit is $u=(0,1)$. 
\end{enumerate}
\end{cor}

\begin{proof} The matrix $\left(\begin{matrix}
 a&b\\
 c&d
\end{matrix}\right)$ fixes $u$ if and only if $b=0$ and $d=1$. Then $a=\pm 1$ and $c$ is arbitrary. Then the conditions (ii) and (iii) of \thmref{main} are equivalent to (i).
\end{proof}

\begin{rem}
 This corollary gives the classification of the irrational noncommutative solenoids of \cite{lp:ncsol} up to isomorphism. This classification appears in Theorem 4.2 of \cite{lp:ncsol}, which considers also rational noncommutative solenoids. The proofs are different.
 \end{rem}

\vskip1cm

\section{ The twisted group $(\Z[1/p]^2,\sigma)$ and the commutation relations.}\label{sectionNCS}

We first consider a general framework: $H$ and $K$ are locally compact abelian groups and $\varphi:H\to\hat K$ is a continuous homomorphism. The groups $H$ and $K$ are noted additively. Their duals $\hat H$ and $\hat K$ are noted multiplicatively. Associated with $\varphi$ we have:
\begin{enumerate}
 \item a 2-cocycle $\sigma\in Z^2(K\times H,\t)$ defined by $\sigma(x,y)=\lan\varphi(x_2),y_1\ran$, where the subscripts $1$ and $2$ stand for the first and second components;
 \item a translation action of $H$ on $\hat K$ defined by $\chi.h=\chi\varphi(h)$.
\end{enumerate}

\begin{prop}\label{CR} The relation $W(k,h)=U(k)V(h)$ gives a bijective correspondence between
\begin{enumerate}
\item pairs $(U,V)$ of unitary representations $U$ of $K$ and $V$ of $H$ satisfying 
\begin{equation}
 U(k)V(h)=<\varphi(h),k>V(h)U(k)
\end{equation}
\item unitary $\sigma$-representations $W$ of $K\times H$ i.e. such that:
\begin{equation}
W(x+y)=\sigma(x,y)W(x)W(y)
\end{equation}
\end{enumerate}
\end{prop}

\begin{proof}
 Suppose $W$ satisfies (ii). Define $U(k)=W(k,0)$ and $V(h)=W(0,h)$. Then
 
\begin{eqnarray*}
V(h)U(k)&=&W(0,h)W(k,0)\cr
&=&\sigma((0,h),(k,0))^{-1}W(k,h)\cr
&=&<\varphi(h),k>^{-1}W(k,h)\cr
U(k)V(h)&=&W(k,0)W(0,h)\cr
&=&\sigma((k,0),(0,h))^{-1}W(k,h)\cr
&=&W(k,h)
\end{eqnarray*}
This gives the relation of (i).\\

Conversely, suppose that $(U,V)$ is as in (i). Define $W(k,h)=U(k)V(h)$. Let us check that $W$ satisfies the relation in (ii). If both $x$ and $y$ belong to the same subgroup $H$ or $K$, $\sigma(x,y)=1$ and the relation is satisfied. Suppose $x=(k,0)$ belongs to $K$ and $y=(0,h)$ belongs to $H$. Then, writing $x=(x_1,x_2)$ and $y=(y_1,y_2)$, we have $x+y=(x_1+y_1,x_2+y_2)$, hence
\begin{eqnarray*}
W(x+y)&=&U(x_1+y_1)V(x_2+y_2)\cr
&=&U(x_1)U(y_1)V(x_2)V(y_2)\cr
&=&U(x_1)<\varphi(x_2),y_1>V(x_2)U(y_1)V(y_2)\cr
&=&<\varphi(x_2),y_1>W(x)W(y)\cr
&=&\sigma(x,y)W(x)W(y)
\end{eqnarray*}

\end{proof}

\begin{thm}\label{isomorphisms} The following C*-algebras are naturally isomorphic
\begin{enumerate}
\item $C^*(K\times H,\sigma)$;
\item $C^*(\hat K\rtimes_\varphi H)$;
\item $C^*(K\ltimes_{\varphi^*}\hat H)$, where $\varphi^*:K\to\hat H$ is the complex conjugate of the transpose of $\varphi$.
 \end{enumerate}
\end{thm}

\begin{proof}
 This is a particular case of \cite[Theorem 7.8]{ren:atge}, to which we refer the reader for further details, but we may as well sketch the proof and give the relevant formulas. Given $f\in C_c(G,\sigma)$, we define its partial Fourier transform. This is the function $\hat f$ on $\hat K\times H$ defined by
$$\hat f(\chi, h)=\int_K f(k,h) \langle\chi,k\rangle dk.$$
It is a continuous function. Its support in the first variable $\chi$ is not compact but its $I$-norm
is bounded. We recall that the convolution product and the involution are defined on the completion $L^I(\hat K\rtimes_\varphi H)$ by the same formulas as on $C_c(\hat K\rtimes_\varphi H)$. Given $f,g\in C_c(G)$, we have
\begin{eqnarray*}
f*g(k,h)&=&\int_G f(k+k',h+h')g(-k',-h')\sigma((k+k',h+h'), (-k',-h')) dk'dh'\cr
&=&\int_{K\times H} f(k+k',h+h')g(-k',-h') \langle\varphi(h+h'),-k'\rangle dk'dh'\cr
\widehat{f*g}(\chi,h)&=&\int f(k+k',h+h')g(-k',-h')\langle\varphi(h+h'),-k'\rangle\langle \chi,k\rangle dk' dh' dk\cr
\hat f*\hat g(\chi,h)&=&\int_H\hat f(\chi,h+h')\hat g(\chi\varphi(h+h'),-h')dh'\cr
&=&\int_H\int_K f(k,h+h') \langle\chi,k\rangle dk\int_K g(k',-h') \langle\chi\varphi(h+h'),k'\rangle dk' dh'\cr
&=&\int f(k,h+h')  g(k',-h')  \langle\chi,k\rangle\langle\chi\varphi(h+h'),k'\rangle dk dk' dh'\cr
&=&\int f(k,h+h')  g(k',-h')  \langle\chi,k\rangle\langle\chi,k'\rangle\langle\varphi(h+h'),k'\rangle  dk dk' dh'\cr
&=&\int f(k,h+h')  g(k',-h')  \langle\chi,k+k'\rangle\langle\varphi(h+h'),k'\rangle  dk dk' dh'\cr
&=&\int f(k,h+h')  g(-k',-h')  \langle\chi,k-k'\rangle\langle\varphi(h+h'),-k'\rangle  dk dk' dh'\cr
&=&\int f(k+k',h+h')  g(-k',-h')  \langle\chi,k\rangle\langle\varphi(h+h'),-k'\rangle  dk dk' dh'\cr
\end{eqnarray*}
We also have
\begin{eqnarray*}
f^*(k,h)&=&\overline{f(-k,-h)}\overline{\sigma((k,h),(-k,-h)}\cr
&=&\overline{f(-k,-h)}\overline{\langle\varphi(h),-k\rangle}\cr
\widehat{f^*}(\chi,h)&=&\int f^*(k,h)\langle\chi,k\rangle dk\cr
&=&\int \overline{f(-k,-h)}\langle \chi\varphi(h),k\rangle dk\cr
(\hat f)^*(\chi,h)&=&\overline{\hat f (\chi\varphi(h),-h)}\cr
&=&\overline{\int f (k,-h)\langle\chi\varphi(h),k\rangle dk}\cr
&=&\overline{\int f (-k,-h)\langle\chi\varphi(h),-k\rangle dk}\cr
&=&\int \overline{ f (-k,-h)}\langle\chi\varphi(h),k\rangle dk\cr
\end{eqnarray*}
The representations of $C^*(\hat K\rtimes_\varphi H)$ are given by a pair $(M,V)$, where
$M$ is a representation of $C_0(\hat K)$, $V$ is a unitary representation of $H$ on the same Hilbert space which satisfy the covariance relation
$$(CR)\quad V(h)^{-1}M(\varphi)V(h)=M(\varphi\circ\alpha_h)$$ 
where $\varphi\in C_0(\hat K)$ and $\alpha_h(\chi)=\chi\varphi(h)$. Via the Gelfand transform,
$M$ corresponds to a unitary representation $U$ of $K$. The correspondence is given by
\[M(\hat f)=\int f(k)U(k)dk\]
where $f\in C_c(K)$ and $\hat f(\chi)=\int f(k)\langle\chi,k\rangle dk$. We show now that $(M,V)$ satisfies (3) iff $(U,V)$ satisfies (1). Suppose (1). Then
\begin{eqnarray*}
V(h)^{-1}M(\hat f)V(h)&=&\int f(k) V(h)^{-1}U(k)V(h) dk\cr
&=&\int f(k) \langle \varphi(h),k\rangle U(k) dk\cr
\end{eqnarray*}
Since the Gelfand transform of $g(k)=f(k) \langle \varphi(h),k\rangle$ is $\hat f\circ\alpha_k$:
\begin{eqnarray*}
\hat g(\chi)&=&\int f(k) \langle \varphi(h),k\rangle\langle\chi,k\rangle dk\cr
&=&\int f(k) \langle\chi \varphi(h),k\rangle dk\cr
&=&\hat f(\chi \varphi(h))\cr
&=&\hat f\circ\alpha_h(\chi)\cr
\end{eqnarray*}
we have (3) for $\varphi=\hat f$. By density, this is true for all $\varphi\in C_0(\hat K)$. Let us prove the converse. We fix $h\in H$. The relation
$$\int f(k)V(h)^{-1}U(k)V(h) dk=\int f(k) \langle \varphi(h),k\rangle U(k)dk$$
for all $f\in C_c(K)$ and the continuity of the functions $k\mapsto V(h)^{-1}U(k)V(h)$ and $k\mapsto  \langle \varphi(h),k\rangle U(k)$ implies the equality
\[V(h)^{-1}U(k)V(h) =  \langle \varphi(h),k\rangle U(k).\]

\end{proof}

The fundamental commutation relation of the noncommutative torus is
\[UV=\zeta\, VU,\]
where $\zeta=e^{2i\pi\alpha}$ and $\alpha\in\R/Z$. Here $H=K=\Z$ and $\varphi:\Z\to\t$ is determined by $\zeta=\varphi(1)$. The representations of this relation are the representations of the twisted group C*-algebra $C^*(\Z^2,\sigma)$, where $\sigma(x,y)=\zeta^{x_1y_2}$. Recall that what matters with projective representations is the antisymmetric bicharacter $\sigma^{(2)}(x,y)=\sigma(x,y)\sigma(y,x)^{-1}=\zeta^{\det(x,y)}$, where $\det(x,y)={x_1y_2-x_2y_1}$.
\vskip 2mm
\begin{lem}\label{generator}
 A homomorphism $\varphi: Z[1/p]\to{\mathcal S}_p$ is completely determined by the image $\zeta=\varphi(1)$ of 1 by $\varphi$. If $\zeta=(\zeta_n)$, then $\lan\varphi(1/p^k), 1/p^l\ran=\zeta_{k+l}$.
\end{lem}
\begin{proof}
 We first show that $\lan\varphi(1/p^k), 1/p^l\ran$ depends on $k+l$ only: replacing $1/p^k$ by $p^l(1/p^{k+l})$ and $1/p^l$ by $p^k(1/p^{k+l})$, one obtains 
 $$\lan\varphi(1/p^k), 1/p^l\ran=\lan\varphi(1/p^{k+l}), 1/p^{k+l}\ran^{(p^{k+l})}.$$
 Let us call it $\zeta(k+l)$. In particular
 $$\zeta(n)=\lan \varphi(1), 1/p^n\ran=\zeta_n.$$
\end{proof}

If $\varphi(1)=\zeta$, we may write $\varphi(x)=\zeta^x$ for $x\in\Z[1/p]$. If  $(\alpha_0,{\bf a}_0)\in\R\times\Q_p$ is a representative of $\zeta$, then $\zeta^x=\pi(x\alpha_0,x{\bf a}_0)$ and this does not depend on the chosen representative. If $x=k/p^n$, where $n\in\N$ and $k\in\Z$, we have $\zeta^x=(S_n\zeta)^k$ where $S_n$ is the $n$-shift. The image of $\varphi$ is
\[\Z[1/p].\zeta\defequal\{\zeta^x, x\in\Z[1/p]\}.\]

\begin{defn} Let $\zeta\in\sol_p$ be an irrational $p$-angle.
\begin{enumerate}
 \item The noncommutative $\zeta$-solenoid is the C*-algebra $A_\zeta=C^*(\Z[1/p]^2,\sigma)$ where $\sigma(x,y)=\langle \zeta, x_1y_2\rangle$.
 \item The irrational rotation C*-algebra of $p$-angle $\zeta$ is the C*-algebra 
 \[\underline A_\zeta=C^*(\Z[1/p]\ltimes_\zeta\sol_p)\] 
 where the translation action of $\Z[1/p]$ on $\sol_p$ is given by the homomorphism $\varphi:\Z[1/p]\to\sol_p$ such that $\varphi(1)=\zeta$.
\end{enumerate}
\end{defn}

If $\zeta=\pi(\alpha,{\bf a})$, we may also use the subscript $(\alpha,{\bf a})$ instead of $\zeta$. With an abuse of language, we shall also say that the groupoid $\Z[1/p]\ltimes_\zeta\sol_p$ is a noncommutative $p$-solenoid.

\vskip2mm 

The natural generators of $Z[1/p]^2$ are $u_k=(1/p^k,0)$ and $v_l=(0,1/p^l)$. They satisfy $u_k=pu_{k+1}$ and $v_l=pv_{l+1}$.

In an analogous fashion to the noncommutative torus, given an $p$-angle $\zeta=(\zeta_n)\in \sol_p$, written as $\hat J'(\alpha,{\bf a})=(e(\alpha/p^n+\{{\bf a}/p^n\}))$, we consider the commutation relations
\[(*)\quad V_lU_k=\zeta_{k+l}U_kV_l;\quad U_k=U_{k+1}^p,\quad V_l=V_{l+1}^p.\]

As a particular case of \propref{CR}, we obtain

\begin{prop}\label{NCSiso} Let $\zeta\in\sol_p$ be an irrational $p$-angle.
\begin{enumerate}
 \item The universal C*-algebra for ($*$) is the twisted group C*-algebra $C^*(\Z[1/p]^2,\sigma)$ where $\sigma(x,y)=\langle \zeta, x_1y_2\rangle$, where we view $\zeta\in\sol_p$ as a character of $\Z[1/p]$.
 \item The C*-algebras $A_\zeta$ and $\underline A_\zeta$ are isomorphic.
\end{enumerate}
 
 \end{prop}
 Note that $\sigma^{(2)}=\zeta\circ\det$ where $\det:\Z[1/p]^2\to\Z[1/p]$ is given by $\det(x,y)={x_1y_2-x_2y_1}$.

\begin{prop} $A_\zeta=\varinjlim B_n$ where $B_n$ is a noncommutative $\zeta_n$-torus.
\end{prop}

\begin{proof} We let $u_k,v_l$, where $k,l\in\N$ be the canonical generators of $A_\zeta$. Since $v_lu_k=\zeta_{k+l}u_kv_l$, the sub C*-algebra $A_{k,l}$ of $A_\zeta$ generated by $u_k$ and $v_l$ is a noncommutative $\zeta_{k+l}$-torus. The system $(A_{k,l}), (k,l)\in\N\times\N$, where $\N\times\N$ is endowed with the product order is inductive and $\bigcup A_{k,l}$ is dense in $A_\zeta$. This implies $A_\zeta=\varinjlim A_{k,l}$. If we wish, we can consider the subsystem $(A_{k,k}), k\in\N$ which has the same limit.
\end{proof}

We find convenient to write $A_\zeta$ as the inductive limit of
\[A_{0,0}\to A_{1,0}\to A_{1,1}\to\cdots\to A_{k,k-1}\to A_{k,k}\to\cdots \]
which we write as
\[B_0\to B_1\to B_2\to\cdots\to B_{2k-1}\to B_{2k}\to\cdots\]
where $ B_{2k-1}= A_{k,k-1}$ and $B_{2k}=A_{k,k}$. Then $B_n$ is a noncommutative torus of parameter $\zeta_n$.

\begin{prop}
  The C*-algebra $A_\zeta$ has a unique tracial state $\tau$. 
\end{prop}

\begin{proof}
See \cite[Remark 7.14]{ren:atge} or \cite[Proposition 32]{gre:loc}.
\end{proof}

We note that the tracial state $\tau$ of $C^*(\Z[1/p]\ltimes_\zeta\sol_p)$ is given by the normalized Haar measure of the compact group $\sol_p$. Using the definition of $A_\zeta$ as inductive limit of irrational rotation C*-algebras, one obtains easily its K-theory.

\begin{prop}\label{K0} Let $(\alpha,{\bf a})\in(\R\setminus\Q)\times\Q_p$ and $\zeta=\pi(\alpha,{\bf a})$. Then $(K_0(A_\zeta),[1])$ is isomorphic to $(G'(\alpha,{\bf a}),(0,1))$ as an ordered group with unit.
 \end{prop}
 
 \begin{proof} Without loss of generality, we may assume that ${\bf a}\in \Z_p$. Since $A_\zeta=\varinjlim B_n$, $K_0(A_\zeta)=\varinjlim K_0(B_n)$. Moreover, the morphisms of the inductive systems preserve the units. The restriction $\tau_n$ of $\tau$ to $B_n$ is its unique tracial state. We know that $\tau_n$ induces an isomorphism of ordered groups from $K_0(B_n)$ onto $\Z\alpha_n+\Z\subset\R$ sending $1$ to $1$. Let $e$ be a projection of $B_n$ of trace $\alpha_n$ and $[e]_n$ its class in $K_0(B_n)$, We  view $e$ as a projection of $B_{n+1}$ defining $[e]_{n+1}$ in $K_0(B_{n+1})$. Its decomposition in the basis $(\alpha_{n+1},1)$ is given by the equation
\[\alpha_n=p\alpha_{n+1}-a_n.\]
If we wish (but this is not necessary to compute it), we can express $K_0(A_\zeta)=\varinjlim K_0(B_n)$ as the inductive limit of the inductive system of groups:
\[\Z^2\to\Z^2\to\cdots\to\Z^2\to\Z^2\to\cdots\]
where the morphism from level $n$ to level $n+1$ is given by the matrix $\left(\begin{matrix}
p&0\\
-a_n&1
\end{matrix}\right)$. However the trace defines immediately an isomorphism of ordered group preserving the unit from $K_0(A_\zeta)$ onto $G(\alpha,{\bf a})=\bigcup_{n\in\N}\Z\alpha_n+\Z$. To conclude, we recall \propref{isomorphism}.
\end{proof}

Note that the state of $(G'(\alpha,{\bf a}),(0,1))$ induced by the tracial state $\tau$ is given by $[\tau](x,y)=x\alpha+y$ for $(x,y)\in G'(\alpha,{\bf a})$.

\vskip1cm

\section{Construction of projective C*-modules.}

Since the fundamental work of M. Rieffel and P. Green on strong Morita equivalence of C*-algebras, and in particular of transformation group C*-algebras (see \cite{rie:morita} and references therein), the construction of C*-modules $C^*(X)$ over $C^*(G)$, where $G$ is a groupoid and $X$ is a proper $G$-space, is well-known (it is briefly recalled in the Appendix). We assume here that the locally compact spaces and the locally compact groupoids of this section are Hausdorff.\\

We first recall the notion of amplification which gives an important example of groupoid equivalence. It is extensively studied in Section 2 of \cite{tu:proper} and is called ``blow up'' in \cite{wil:toolkit}.

\begin{prop}\cite[Proposition 2.7]{tu:proper}
 Let $H$ be a locally compact groupoid, $X$ a locally compact space, and $\pi: X\to H^{(0)}$ a continuous, open and surjective map. Then
 \[X*H*X\defequal\{(x,h,y)\in X\times H\times X: \pi(x)=r(h), \pi(y)=s(h)\}\]
 has a natural structure of locally compact groupoid with unit space $X$.
\end{prop}

The multiplication law is $(x,h,y)(y,k,z)=(x,hk,z)$. The topology is the subset topology of $X\times H\times X$. We say that $X*H*X$ is the {\em amplification} of $H$ along $\pi$. The space
\[X*H=\{(x,h)\in X\times H: \pi(x)=r(h)\}\]
is a groupoid equivalence between $X*H*X$ and $H$ (\cite[Proposition 2.29]{tu:proper} and \cite [Example 2.37]{wil:toolkit}).

We are interested in the following situation: the locally compact group $\Gamma$ acts continuously on the locally compact space $X$ and $\Gamma_1$ is a closed normal subgroup of $\Gamma$ such that the restriction of the action to $\Gamma_1$ is free and proper. We denote by $\rho:\Gamma\to\Gamma/\Gamma_1$ and $\pi: X\to X/\Gamma_1$ the quotient maps. Then $\Lambda=\Gamma/\Gamma_1$ acts continuously on $X/\Gamma_1$ according to $\rho(g)\pi(x)=\pi(gx)$ for $(g,x)\in G\times X$.

\begin{thm}\label{reduction} Assume that the locally compact group $\Gamma$ acts continuously on the locally compact space $X$ and $\Gamma_1$ is a closed normal subgroup of $\Gamma$ such that the restriction of the action to $\Gamma_1$ is free and proper.  Then the semidirect product groupoid $G=\Gamma\ltimes X$ is isomorphic to the amplification of the semidirect product groupoid $H=\Gamma/\Gamma_1 \ltimes X/\Gamma_1$ along $\pi$.\\
\end{thm}

\begin{proof} We first observe that the amplification $X*H*X$ can be written as
\[X*\Lambda*X=\{(y,\lambda,x)\in X\times\Lambda\times X: \pi(y)=\lambda\pi(x)\}.\]
We consider the map
\begin{eqnarray*}
 f:\Gamma\times X&\to X*\Lambda*X\cr
 (g,x)&\mapsto (gx,\rho(g),x)\cr
\end{eqnarray*}
The map $f$ is surjective: given $(y,\lambda,x)\in X\times\Lambda\times X$ such that $\pi(y)=\lambda\pi(x)$, we pick $g\in\Gamma$ such that $\lambda=\rho(g)$. We then have $\pi(y)=\pi(gx)$. Therefore, there exists $g_1\in\Gamma_1$ such that $y=g_1gx$. Then $(y,\lambda,x)=f(g_1g,x)$.\\
The map is injective: suppose that $(g'y,\rho(g'),y)=(gx,\rho(g),x)$. Then, first $y=x$. The equation $\rho(g')=\rho(g)$ gives $g'=g_1g$ for some $g_1\in\Gamma_1$. The equation $g_1gx=gx$ gives $g_1=e$ since $\Gamma_1$ acts freely on $X$. Hence $g'=g$.\\
We leave to the reader to check that the map $f$ is a groupoid homomorphism.\\
The map $f$ is clearly continuous. To show that $f$ is a homeomorphism, we show that $f$ is proper. Let $A,B$ be compact subsets of $X$ and $K$ be a compact subset of $\Lambda$. There exists a compact subset $L$ of $\Gamma$ such that $K\subset \rho(L)$. Since $\Gamma_1$ acts properly on $X$, the set
\[K_1=\{g_1\in\Gamma_1: g_1^{-1}L^{-1}B\cap A\not=\emptyset\}\]
is compact. Then, $f^{-1}(B\times K\times A)$ is contained in the compact set $LK_1\times A$: suppose indeed that $x\in A$, $gx\in B$ and $\rho(g)\in K$. There exist $g_1\in\Gamma_1$ and $g'\in L$ such that $g=g'g_1$. The condition $g'g_1x\in B$ implies that $g_1x\in L^{-1}B$. Then $x$ belongs to $g_1^{-1}L^{-1}B\cap A$, $g_1$ belongs to $K_1$ and $g=g'g_1$ belongs to $LK_1$.
\end{proof}

\begin{cor}\label{cor} With above notation, the space 
\[X*H=\{(x,\lambda,\pi(y)): x,y\in X, \lambda\in\Lambda,\quad \pi(x)=\lambda \pi(y)\}=X\times\Lambda,\]
endowed with the moment maps  $r(x,\lambda)=x$ and $s(x,\lambda)=\lambda^{-1}\pi(x)$, the left $G$-action $g(x, \lambda)=(gx,\rho(g)\lambda)$ for $g\in\Gamma$ and the right $H$-action $(x,\lambda)\mu=(x,\lambda\mu)$ for $\mu\in\Lambda$ is a   $(G,H)$-equivalence.
 \end{cor}

 In our case, $X$ is a group and the action of $\Gamma$ on $X$ is given by a group homomorphism $j:\Gamma\to X$, namely $g x=j(g)x$. Then, the action of $\Gamma/\Gamma_1$ on $X/\Gamma_1$ is given by the induced homomorphism $j':\Gamma/\Gamma_1\to X/\Gamma_1$.\\
 
\begin{prop}\label{groupoid equivalence} Consider the translation action of $\Gamma=\Z[1/p]^2$ on $X=\R\times\Q_p$ given by the homomorphism 
\[J_{(\alpha,{\bf a})}:\Gamma\to\R\times\Q_p: (x,y)\mapsto (x\alpha+y,x{\bf a}-y),\] 
where $(\alpha,{\bf a})\in(\R\setminus\Q)\times\Q_p$ and let $G$ be semidirect product $\Gamma\ltimes_{(\alpha,{\bf a})}X$. Then
\begin{enumerate}
 \item For every $(a,b)\in\Z[1/p]^2$ such that $a{\bf a}-b\not=0$, the action of $\Gamma_1=\Z[1/p](a,b)$ on $X=\R\times\Q_p$ is free and proper. Therefore $X\times(\Gamma/\Gamma_1)$ is a free and proper left $G$-space which gives an equivalence between $G$ and $(\Gamma/\Gamma_1)\ltimes (X/\Gamma_1)$.
 \item If moreover there exist $c,d\in\Z[1/p]$ such that $ad-bc=1$, $(\Gamma/\Gamma_1)\ltimes (X/\Gamma_1)$ is isomorphic to the noncommutative $p$-solenoid $\Z[1/p]\ltimes_{(\beta,{\bf b})}\sol_p$ where
\[(\beta, {\bf b})=\left(\begin{matrix}
 c&d\\
 a&b
\end{matrix}\right)(\alpha, {\bf a}).\]
\end{enumerate}
\end{prop}

\begin{proof}
 For (i), we show that $\Gamma_1=\Z[1/p](a,b)$ satisfies the assumption of \thmref{reduction}. Multiplication by $(a\alpha+b,-a{\bf a}+b)$ is an automorphism of the locally compact additive group $\R\times\Q_p$ which carries the subgroup $J_{(\alpha,{\bf a})}(\Z[1/p](0,1))=\Z[1/p](1,-1)$ onto the subgroup $J_{(\alpha,{\bf a})}(\Z[1/p](a,b))=\Z[1/p](a\alpha+b,a{\bf a}-b)$. For further reference, we note that the module of this automorphism is $|a\alpha+b|\,|-a{\bf a}+b|_p$. We have seen that $\Z[1/p](1,-1)$ is a discrete subgroup of $\R\times\Q_p$. We deduce that $J_{(\alpha,{\bf a})}(\Z[1/p](a,b)=\Z[1/p](a\alpha+b,a{\bf a}-b)$, image of $\Z[1/p](1,-1)$ by this automorphism, is a discrete subgroup of $\R\times\Q_p$. Therefore, it acts freely and properly on $\R\times\Q_p$.\\
 For (ii), we identify $\Gamma/\Gamma_1$ to $\Z[1/p]$ through the homomorphism $(x,y)\mapsto -bx+ay$ and $X/\Gamma_1$ to $\sol_p$ through the homomorphism $(t,x)\mapsto \pi(t/r,x/\rho)$, where $\pi:\R\times\Q_p\to\sol_p$ is the usual quotient map, $r=a\alpha+b$ and $\rho=a{\bf a}-b$. Then the induced homomorphism $j':\Gamma/\Gamma_1\to X/\Gamma_1$ becomes a homomorphism $j'':\Z[1/p]\to\sol_p$. According to \lemref{generator}, such a homomorphism is determined by $j''(1)$.
We deduce from the relation $ad-bc=1$ that $j''(1)=\pi(\beta,{\bf b})$, where $\beta=\dfrac{c\alpha+d}{a\alpha+b}$ and ${\bf b}=\dfrac{c{\bf a}-d}{a{\bf a}-b}$.
\end{proof}

In particular, when $(a,b,c,d)=(0,1,1,0)$, one obtains the noncommutative $p$-solenoid $\Z[1/p]\ltimes_{(\alpha,{\bf a})}\sol_p$.\\

 The groupoid equivalence $X\times(\Gamma/\Gamma_1)$ between $G=\Gamma\ltimes_{\alpha,{\bf a}}X$ and $H=\Gamma/\Gamma_1\ltimes X/\Gamma_1$ implements a Morita equivalence of their C*-algebras and an isomorphism of their $K$-theory groups. The identification of the $K_0$-groups is made easy through the trace. We proceed as follows. We first recall that a Radon measure $\mu$ on the unit space $G^{(0)}$ of a locally compact groupoid with Haar system $(G,\lambda)$ which is invariant in the sense that $\mu\circ\lambda=(\mu\circ\lambda)^{-1}$
defines a trace $\tau_\mu$ on the full C*-algebra $C^*(G,\lambda)$ such that $\tau_\mu(f)=\int_{G^{(0)}}f_{|G^{(0)}} d\mu$ for $f\in C_c(G)$ and a transverse measure $\Lambda_\mu$ on the category of principal $G$-spaces (see \cite[Appendix A.1, Section {\bf a}]{adr:amenable}).

\begin{prop}\label{dimension} Assume that $\Gamma$ is a discrete group acting continuously on the locally compact space $X$ and that $\Gamma_1$ is a normal subgroup whose action on $X$ is free and proper. Let $\mu$ be a Radon measure on $X$ invariant with respect to the groupoid $G=\Gamma\ltimes X$.  If $X/\Gamma_1$ is compact, then the $\tau_\mu$-dimension of the projective C*-module $C^*(X\times(\Gamma/\Gamma_1))$ is $\Lambda_\mu(X/\Gamma_1)$.
 \end{prop}

\begin{proof}
 This is a particular case of \cite[Proposition 5.5.6]{ren:rio} which is reproduced in the Appendix as \propref{volume}.
\end{proof}

In our situation, where  $\Gamma=\Z[1/p]^2$ acts on $X=\R\times\Q_p$ by translation $J_{(\alpha,{\bf a})}$, where $(\alpha,{\bf a})\in(\R\setminus\Q)\times\Q_p$, we consider the measure
$\mu=\lambda\times\lambda_p$, where $\lambda$ is the Lebesgue measure on $\R$ and $\lambda_p$ is the normalized Haar measure of $(\Q_p,+)$.

\begin{cor} Let $\Gamma, X, J_{(\alpha,{\bf a})}, \mu$ be as above. Let  $\Gamma_1=\Z[1/p](a,b)$ where $(a,b)\in\Z[1/p]^2$ satisfies $a{\bf a}-b\not=0$. Then the $\tau_\mu$-dimension of the projective C*-module $C^*(X\times (\Gamma/\Gamma_1))$ is $|a\alpha+b||a{\bf a}-b|_p$.
\end{cor}

\begin{proof} We have seen that $\Z[1/p](1,-1)$ is a discrete subgroup of $\R\times\Q_p$. The quotient is $\sol_p$ which has volume 1. The subgroup $J_{(\alpha,{\bf a})}(\Z[1/p](a,b))=\Z[1/p](a\alpha+b,a{\bf a}-b)$ is the image of $\Z[1/p](1,-1)$ by the automorphism of multiplication by $(a\alpha+b,-a{\bf a}+b)$ which has module $|a\alpha+b|\,|-a{\bf a}+b|_p$. The volume of the quotient is the measure of a fundamental domain. We have seen that $[0,1)\times\Z_p$ is a fundamental domain for $\R\times\Q_p/\Z[1/p](1,-1)$. Therefore, $[0,a\alpha+b)\times (a{\bf a}-b)\Z_p$ is a fundamental domain for $\R\times\Q_p/\Z[1/p](a\alpha+b,a{\bf a}-b)$. Hence, the volume of the quotient is  $|a\alpha+b||a{\bf a}-b|_p$. In conclusion, $Z(a,b)=X\times(\Gamma/\Z[1/p](a,b))$ is a principal $G$-space, where $G=\Gamma\ltimes X$, giving the projective C*-module $C^*(Z(a,b))$ of dimension  $|a\alpha+b|\,|a{\bf a}-b|_p$. 
\end{proof}

We denote by $[C^*(Z(a,b))]$ its class in the K$_0$-group of $C^*(G)$. Thus, we define a map
\[\Phi:\{(a,b)\in\Gamma: a{\bf a}-b\not=0\}\to K_0^+(C^*(G)): (a,b)\mapsto[C^*(Z(a,b))]\]
such that $\tau_*(\Phi(a,b))=|a\alpha+b|\,|a{\bf a}-b|_p$, where $\tau_*$ is the homomorphism induced by the trace $\tau_\mu$. Since, according to \propref{groupoid equivalence}, $C^*(G)$ is Morita equivalent to the C*-algebra of the noncommutative $p$-solenoid $\Z[1/p]\ltimes_{(\alpha,{\bf a})}\sol_p$, their positive cones $K_0^+$ are isomorphic. It has been shown in \propref{K0} that  $K_0^+(C^*(\Z[1/p]\ltimes_{(\alpha,{\bf a})}\sol_p))$ is isomorphic to $G'(\alpha,{\bf a})^+$. The map $\Phi$ has then a natural expression given in the following proposition.

\begin{prop}\label{Psi} Given $(a,b)\in\Gamma$ such that $ a{\bf a}-b\not=0$, $\Psi(a,b)\defequal \epsilon p^{-n}(a,b)$, where $n=v(a{\bf a}-b)$ and $\epsilon={\rm sgn}(a\alpha+b)$, belongs to $G'(\alpha,{\bf a})^+$. The map
 \[\Psi:\{(a,b)\in\Gamma: a{\bf a}-b\not=0\}\to G'(\alpha,{\bf a})^+\]
 satisfies
 \begin{enumerate}
 \item $j_\alpha\circ\Psi=\tau_*\circ\Phi$;
 \item $\Psi(a',b')=\Psi(a,b)$ if and only if there exists an invertible element $\lambda\in\Z[1/p]$ such that $(a',b')=\lambda(a,b)$;
 \item its image is
 \[\{(a,b)\in\Z[1/p]\times\Z[1/p]: v(a{\bf a}-b)=0\quad{\rm and}\quad a\alpha+b\ge 0\}.\]
\end{enumerate}
\end{prop}

\begin{proof} It is immediate to check that $(a',b')=\Psi(a,b)$ satisfies $v(a'{\bf a}-b')=0$ and $a'\alpha+b'=p^{-n}|a\alpha+b|\ge 0$.\\ 
The property (i) is clear:
\begin{eqnarray*}
 j_\alpha\circ\Psi(a,b)&=p^{-n}(a\alpha+b)\cr
 \tau_*\circ\Phi(a,b)&=|a\alpha+b|p^{-n}
\end{eqnarray*}
Assertion (ii):
if $\Psi(a',b')=\Psi(a,b)$, then $(a',b')=\lambda(a,b)$ where, with the notation of the proposition, $\lambda=(\epsilon/\epsilon')p^{n-n'}$. Conversely,  $\Psi(\lambda a,\lambda b)=\Psi(a,b)$ for $\lambda=\epsilon p^k$ where $\epsilon=\pm 1$ and $k\in\Z$.\\
Assertion (iii): we have seen that  $(a',b')=\Psi(a,b)$ satisfies $v(a'{\bf a}-b')=0$ and $a'\alpha+b'=p^{-n}|a\alpha+b|\ge 0$. Every $(a,b)\in \Z[1/p]\times\Z[1/p]$ such that $v(xa{\bf a}-b)=0$ and $a\alpha+b\ge 0$ is in the image since then $(a,b)=\Phi(a.b)$.
\end{proof}

We return to the general framework of a locally compact group $\Gamma$ acting continuously on a locally compact space $X$. We assume now that $\Gamma_1$ and $\Gamma_2$ are two closed normal subgroups of $\Gamma$ which act freely and properly on $X$. We denote by $\rho_i: \Gamma\to\Gamma/\Gamma_i$ and $\pi_i:X\to X/\Gamma_i$ the quotient maps. We also write $\Lambda_i=\Gamma/\Gamma_i$, $H_i=\Lambda_i\ltimes X/\Gamma_i$ and $Z_i=X\times\Lambda_i$, which is a $(G,H_i)$-equivalence. By composition, one obtains a $(H_1,H_2)$-equivalence. We want to understand better this equivalence as a topological space. Instead of introducing the inverse equivalence $\tilde Z_1$ and considering $\tilde Z_1\circ Z_2$, we may as well describe it as $Z= (Z_1*Z_2)/G$, where
$Z_1*Z_2=\{(z_1,z_2): r(z_1)=r(z_2)\}$ and $G$ acts by $\gamma(z_1,z_2)=(\gamma z_1,\gamma z_2)$. Then the left action of $H_1$ and the right action of $H_2$ are given by
\[\eta_1[z_1,z_2]\mu_2=[z_1\eta_1^{-1}, z_2\mu_2], \quad \mu_1\in \Lambda_1,\,\mu_2\in \Lambda_2.\]
Explicitly, we can write this $(H_1,H_2)$-equivalence as
$Z=(\Lambda_1\times X\times\Lambda_2)/\Gamma$,
where 
$g(\lambda_1,x, \lambda_2)=(\lambda_1\rho_1(g)^{-1},gx,\rho_2(g)\lambda_2)$ for $g\in\Gamma$ and $\lambda_i\in \Lambda_i$.
The moment maps $r:Z\to X/\Gamma_1$ and $s:Z\to X/\Gamma_2$ are given respectively by
$r[\lambda_1,x,\lambda_2]=\lambda_1\pi_1(x)$ and $s[\lambda_1,x,\lambda_2]=\lambda_2^{-1}\pi_2(x)$. The left and right actions are given by\\
$\mu_1[\lambda_1,x,\lambda_2]\mu_2=[\mu_1\lambda_1,x,\lambda_2\mu_2]$ for $\mu_i\in \Lambda_i$.\\


It is interesting to give other presentations of this $(H_1,H_2)$-equivalence. We define

$Z'=(\Gamma\times X)/(\Gamma_1\times\Gamma_2)$ where $\Gamma_1\times\Gamma_2$ acts on $\Gamma\times X$ by 
$(g_1,g_2)(g,x)=(g_1g g_2^{-1}, g_2x)$.
The topology is the quotient topology.
The moment maps $r:Z'\to X/\Gamma_1$ and $s:Z'\to X/\Gamma_2$ are given respectively
by $r[g,x]=\pi_1(gx)$ and $s[g,x]=\pi_2(x)$.
The left and right actions are given by
$\rho_1(g')[g,x]=[g'g,x]$ and $[g,x]\rho_2(g')=[gg',g'^{-1}x]$.\\

We also define $Z''=\{(\pi_1(y),\pi_2(x)): x,y\in X\quad{\rm and}\quad \exists g\in\Gamma: y=gx\}$.
The topology is the final topology with respect to the map from $\Gamma\times X$ to $Z''$ sending $(g,x)$ to $(\pi_1(gx),\pi_2(x))$. The moment maps $r:Z''\to X/\Gamma_1$ and $s:Z''\to X/\Gamma_2$ are respectively
$r(x_1,x_2)=x_1$ and $s(x_1,x_2)=x_2$.
The left and right actions are given by
$\rho_1(g)(x_1,x_2)=(\rho_1(g)x_1,x_2)$ and $(x_1,x_2)\rho_2(g)=(x_1,\rho_2(g)^{-1}x_2)$.

\begin{prop}  Let $\Gamma, X, \Gamma_1,\Gamma_2$ be as above. We assume that $\Lambda_1=\Gamma/\Gamma_1$ and $\Lambda_2=\Gamma/\Gamma_2$ are discrete.
\begin{enumerate}
\item $Z$ and $Z'$ are isomorphic $(H_1, H_2)$-equivalences.
\item If the action of $\Gamma$ on $X$ is free, these $(H_1, H_2)$-equivalences are also isomorphic to $Z''$.
\end{enumerate}
 \end{prop}

\begin{proof}
 (i) We show that the map
$\varphi:\Gamma\times X\to\Lambda_1\times X\times\Lambda_2$ sending
$(g,x)$ to $(\rho_1(g),x,e)$ induces a bijection from $Z'$ onto $Z$. Since $(g_1,g_2)(g,x)=(g_1g g_2^{-1}, g_2x)$ is sent to
 \[(\rho_1(g_1g g_2^{-1}),g_2x,e)=(\rho_1(g)\rho_1(g_2)^{-1},g_2x,\rho_2(g_2)e)=g_2(\rho_1(g),x,e),\] the map passes to the quotients giving
\begin{eqnarray*}
\dot\varphi:(\Gamma\times X)/(\Gamma_1\times\Gamma_2)&\to& (\Lambda_1\times X\times\Lambda_2)/\Gamma\cr
[g,x]&\mapsto& [\rho_1(g),x,e]\cr
\end{eqnarray*}
Let us show that $\dot\varphi$ is surjective: given $(\lambda_1,y,\lambda_2)$, we write $\lambda_i=\rho_i(h_i)$ where $h_i\in\Gamma$. Then
 \[(\lambda_1, y,\lambda_2)= (\rho_1(h_1),y,\rho_2(h_2))=h_2 (\rho_1(h_1h_2),h_2^{-1}y,e).\]
Therefore, $[\lambda_1, y,\lambda_2]=\dot\varphi[g,x]$, where $g=h_1h_2$ and $x=h_2^{-1}y$.\\
 Let us show that $\dot\varphi$ is injective. Suppose that $[\rho_1(g),x,e]=[\rho_1(h),y,e]$. There exists $g_2\in\Gamma$ such that 
\begin{eqnarray*}
 (\rho_1(h),y,e)&=&g_2(\rho_1(g),x,e)\cr
 &=&(\rho_1(g)\rho_1(g_2)^{-1},g_2x,\rho_2(g_2))\cr
 &=&(\rho_1(gg_2^{-1}),g_2x,\rho_2(g_2))\cr
\end{eqnarray*}
This implies that $g_2\in \Gamma_2$, $y=g_2x$ and $h=g_1gg_2^{-1}$ for some $g_1\in\Gamma_1$. Thus, 
 $$(h,y)=(g_1gg_2^{-1},g_2x)=(g_1,g_2)(g,x)$$
 and $[h,y]=[g,x]$. Since $\varphi$ is continuous and open (because we assume $\Lambda_i$ discrete), so is $\dot\varphi$. Therefore, it is a homeomorphism. It is straightforward to check that this homeomorphism respects the structure of correspondence.\\

  (ii) Two elements $(g,x)$ and $(h,y)$ in $\Gamma\times X$ have the same image in $Z''$ if and only if there exist $(g_1,g_2)\in \Gamma_1\times \Gamma_2$ such that $y=g_2x$ and $hy=g_1gx$. If the action of $\Gamma$ on $X$ is free, this is equivalent to $y=g_2x$ and $h=g_1gg_2^{-1}$. In other words, $Z''$ and $Z'$ are the quotient spaces of the same equivalence relation on $\Gamma\times X$. Therefore, the map from $Z'$ to $Z''$ sending $[g,x]$ to $(\pi_1(gx),\pi_2(x))$ is a homeomorphism. It is straightforward to check that is an isomorphism of correspondences from $Z$ onto $Z''$.\\

\end{proof}

\begin{prop}\label{index}
 Assume that $\Gamma$ is a discrete abelian group acting freely and continuously on the locally compact space $X$. Assume that $\Gamma_1$ and $\Gamma_2$ are subgroups acting properly on $X$ and that $\Gamma_1\cap\Gamma_2=\{e\}$. Then $Z'$ is the disjoint union over $\Gamma/\Gamma_1\Gamma_2$ of open subsets homeomorphic to $X$.
\end{prop}

\begin{proof}
 We denote by $[g,x]$ the class of $(g,x)$ in $Z'=(\Gamma\times X)/(\Gamma_1\times\Gamma_2)$. For a fixed $g\in\Gamma$, the map $x\mapsto [g,x]$ from $X$ to $Z'$ is injective: if [g,y]=[g,x], there exists $(g_1,g_2)\in\Gamma_1\times\Gamma_2$ such that $(g,y)=(g_1gg_2^{-1},g_2x)$. Since implies $g_1=g_2=e$ and $y=x$. Since this map is open, its image is an open subset $[g,X]$ of $Z'$ and $Z'$ is the union over $\Gamma$ of these open subsets. To conclude, we show that, if $g$ and $h$ are congruent modulo $\Gamma_1\Gamma_2$, then $[g,X]=[h,X]$ and, if not, $[g,X]$ and $[h,X]$ are disjoint. Indeed, if $h=g_1gg_2^{-1}$ with $(g_1,g_2)\in\Gamma_1\times\Gamma_2$, then for all $y\in X$, $[h,y]=[g,g_2^{-1}y]$. On the other hand, if $[g,x]=[h,y]$, there exists $(g_1,g_2)\in\Gamma_1\times\Gamma_2$ such that $h=g_1hg_2^{-1}$.
 \end{proof}
 
 We apply \propref{index} to our example: $\Gamma=\Z[1/p]\times\Z[1/p]$ acts on $X=\R\times\Q_p$  by
\[J_{(\alpha,{\bf a})}:\Gamma\to\R\times\Q_p: (x,y)\mapsto (x\alpha+y,x{\bf a}-y),\] 
where $(\alpha,{\bf a})\in (0,1)\setminus\Q\times\Z_p$. We let $\Gamma_i=\Z[1/p](a_i,b_i)$, $i=1,2$, where $(a_i,b_i)\in\Gamma$. If $(a_1,b_1)$ and $(a_2,b_2)$ are linearly independent, $\Gamma_1\cap\Gamma_2=\{(0,0)\}$ and, by the structure theorem of finitely generated modules over a principal ideal domain, $\Gamma_1+\Gamma_2$ has finite index $|n|$, where
\[\det\left(\begin{matrix}
a_1&a_2\\
b_1&b_2
\end{matrix}\right)=\dfrac{n}{p^k},\quad n\wedge p=1.\]
Therefore, the $(H_1,H_2)$-equivalence $Z$, where $H_i=\Gamma/\Gamma_i\ltimes X/\Gamma_i$ is a disjoint union of $|n|$ copies of $\R\times\Q_p$.\\

If $\Gamma_1=\Z[1/p](a_1,b_1)$ where $a_1,b_1\in\Z[1/p]$ are relatively prime in the sense that there exist $a_2, b_2\in\Z[1/p]$ such that $a_1b_2-a_2b_1=1$, then $H=\Gamma/\Gamma_1\ltimes X/\Gamma_1$ is isomorphic to an irrational noncommutative $p$=solenoid $\Z[1/p]\ltimes\sol_p$. Explicitly, let us first choose $\Gamma_1=\Z[1/p](0,1)$. The induced homomorphism $j':\Gamma/\Gamma_1\to X/\Gamma_1$ of \remref{translation} can be identified to the homomorphism $\varphi:\Z[1/p]\to\sol_p$ such that $\varphi(1)=\pi(\alpha,{\bf a})$ and $\Gamma/\Gamma_1\ltimes X/\Gamma_1$ can be identified to $\Z[1/p]\ltimes_{(\alpha,{\bf a})}\sol_p$. Let $\Gamma_2=\Z[1/p](a,b)$ where $a,b\in\Z[1/p]$ are relatively prime with $ad-bc=1$. Then the automorphism $\left(\begin{matrix}
 c&a\\
 d&b
\end{matrix}\right)$ carries $\Gamma_1$ to $\Gamma_2$. The commutative diagram
 \[\begin{CD}
\Z[1/p]\times\Z[1/p]@>{J_{(\alpha,{\bf a})}}>{}>\R\times Q_p\\
@V{\left(\begin{matrix}
 c&a\\
 d&b
\end{matrix}\right)}V{}V  @V{(r,\rho)}V{}V\\
\Z[1/p]\times\Z[1/p] @>{J_{(\beta,{\bf b})}}>{}>\R\times Q_p
\end{CD}\]
shows that
\[(\alpha, {\bf a})=\left(\begin{matrix}
c&d\\
 a&b
\end{matrix}\right)(\beta, {\bf b}).\]
Thus, $\Gamma/\Gamma_2\ltimes X/\Gamma_2$ can be identified to $\Z[1/p]\ltimes_{(\beta,{\bf b})}\sol_p$. The $(H_1,H_2)$-equivalence $Z$ is a disjoint union of $|n|$ copies of $\R\times\Q_p$, where $a=\dfrac{n}{p^k}$ with $n\wedge p=1$.

\section{Noncommutative solenoids and groupoids, Morita equivalence classification}

In this last section we will summarize our Morita equivalence classification results in a single statement.

Let $(\alpha,{\bf a})\in (0,1)\times\Z_p$ irrational. There are at least two groupoids that have already appeared in the previous section:

\begin{enumerate}
\item The (big) solenoidal groupoid
\begin{equation}
\mathcal{G}_{(\alpha,{\bf a})}:= \Z[1/p]^2 \ltimes_{J_{(\alpha,{\bf a})}} (\R\times \Q_p)
\end{equation}
associated to the action by the group morphism 
$$J_{(\alpha,{\bf a})}:\Z[1/p]^2\to \R\times \Q_p$$
given by $J_{(\alpha,{\bf a})}(x,y):=(x\alpha+y,x{\bf a}-y)$.
\item The (small) solenoidal groupoid
\begin{equation}
\sol_{(\alpha,{\bf a})}:= \Z[1/p] \ltimes_{(\alpha,{\bf a})} {\mathcal S}_p
\end{equation}
associated to the action 
$$m\cdot z:=\pi(m\alpha,m{\bf a})\cdot z.$$
\end{enumerate}
As discussed in the previous section (specially proposition \ref{groupoid equivalence}), these two groupoids are equivalent, more precisely the big one is an amplification of the small one. The groupoid $C^*$-algebra of the small one is isomorphic to the noncommutative solenoid $A_{\pi((\alpha,{\bf a}))}$, proposition \ref{NCSiso}.

We can state our final theorem, it resumes the results of the previous sections in a sinlge statement.

\begin{thm}\label{thm1grpdsbis} Let $(\alpha,{\bf a}),(\beta,{\bf b})\in (0,1)\times\Z_p$ irrationals. Then \tfae
\begin{enumerate}
\item The noncommutative solenoids $A_{\pi(\alpha,{\bf a})}$ and $A_{\pi(\beta,{\bf b})}$ are strong Morita equivalent $C^*$-algebras.
 \item The dimension groups $G[\alpha,{\bf a}]$ and $G[\beta,{\bf b}]$ are isomorphic as ordered groups.
 \item  There exist $k,l\in\N$ such that $\sol_k(\alpha, {\bf a})$ and $\sol_l(\beta, {\bf b})$ are congruent modulo ${\rm GL}(2,\Z)$.
 \item $(\alpha, {\bf a})$ and $(\beta, {\bf b})$ are congruent modulo ${\rm GL}(2,\Z[1/p])$.
 \item The (small) solenoidal groupoids $\sol_{(\alpha,{\bf a})}$ and $\sol_{(\beta,{\bf b})}$ are Morita equivalent.
\end{enumerate}
 \end{thm}
 
 \begin{proof}
 $(i)\implies (ii):$ This is obtained in section \ref{sectionNCS}, more explicitly, it is an immediate consequence of proposition \ref{K0} since a strong Morita equivalence between $A_{\pi((\alpha,{\bf a}))}$ and $A_{\pi((\beta,{\bf b}))}$ induces an isomorphism of their $K-$theory groups as ordered groups.
 
 $(ii)\iff (iii):$ It is the statement of \corref{thm1} and \remref{shift}.
 
 $(iii)\implies (iv):$ It is a direct consequence of proposition \ref{propcongruent}.$(iii)$.
 
 $(iv)\implies (v):$ This point is a corollary of the last proposition of the previous section, let us repeat the argument in terms of groupoids. Let us assume that $(\alpha, {\bf a})$ and $(\beta, {\bf b})$ are congruent through the matrix $A= \left(\begin{matrix}
 a&c\\
 b&d
\end{matrix}\right)\in{\rm GL}(2,\Z[1/p])$, that is 
\[(\alpha, {\bf a})=A \cdot (\beta, {\bf b}),\]
as in definition \ref{defcongruent}. Consider the following diagram
\begin{equation}
\xymatrix{
\Z[1/p]^2\ar[d]_-{^tA}\ar[rr]^-{J_{(\alpha,{\bf a})}}&& \R\times \Q_p\ar[d]^-{\cdot (b\beta + d,-b{\bf b}+d)}\\
\Z[1/p]^2\ar[rr]_-{J_{(\beta,{\bf b})}}&& \R\times \Q_p
}
\end{equation}
where $\cdot (b\beta + d,-b{\bf b}+d)$ stands for the multiplication by the invertible element $(b\beta + d,-b{\bf b}+d)$ of $\R\times \Q_p$, invertible since $(\beta,{\bf b})$ is irrational. It is straightforward to verify that the diagram above is commutative and since both vertical maps are isomorphisms we deduce an induced groupoid isomorphism between the associated big solenoidal groupoids
\begin{equation}
\mathcal{G}_{(\alpha,{\bf a})}\cong \mathcal{G}_{(\beta,{\bf b})}
\end{equation}
and hence an equivalence of groupoids between the associated small solenoidal groupoids
\begin{equation}
\sol_{(\alpha,{\bf a})}\sim \mathcal{G}_{(\alpha,{\bf a})}\cong \mathcal{G}_{(\beta,{\bf b})}\sim \sol_{(\beta,{\bf b})}.
\end{equation}
$(v)\implies (i):$ If $\sol_{(\alpha,{\bf a})}\sim \sol_{(\beta,{\bf b})}$ are Morita equivalent groupoids, then they induce a strong Morita equivalence of their associated $C^*$-algebras.
 \end{proof}
\vskip 1cm
\appendix
\section{Explicit construction of projections}

For the convenience of the reader, we recall the definition of the module operation and the inner product of the C*-module $C^*(Z)$ of a proper locally compact right $G$-space $Z$, where $G$ is a locally compact groupoid. For the sake of simplicity, we assume that $Z$ and $G$ are Hausdorff (the non-Hausdorff case is studied in \cite{tu:proper}). The construction requires that $G$ is equipped with a Haar system $(\lambda^x)_{x\in G^{(0)}}$ and $Z$ is equipped with a $G$-equivariant system of measures $(\alpha_x)_{x\in G^{(0)}}$ along the fibres of the moment map $s_Z:Z\to G^{(0)}$. Such a pair $(Z,\alpha)$ is called a {\it transverse function} in \cite[Appendix A.1]{adr:amenable}). The details can be found in \cite{mrw:equivalence, ren:rep, tu:proper}. For $h,k\in C_c(Z)$ and $f\in C_c(G)$, we define
\[h*f(z)=\int_G h(z\gamma)f(\gamma^{-1})d\lambda^{s_Z(z)}(\gamma),\qquad \<h,k\>(\gamma)=\int_Z \overline{h(z\gamma^{-1})}k(z)d\alpha_{r(\gamma)}(z).\]
The next proposition is a mild generalisation (from ``$G$ \'etale'' to ``$G$ has a Haar system'') of \cite[Corollary 5.5.5]{ren:rio}. It appears in the framework of foliation algebras in \cite[Section 7]{con:survey} and \cite[II.8.$\beta$]{con:ncg} but also implicitly in \cite{rie:irrational}.

\begin{prop}\label{projection} Let $(G,\lambda)$ be a locally compact Hausdorff groupoid with Haar system and let $Z$ be a free and proper locally compact Hausdorff right $G$-space. Assume that the moment map $s_Z:Z\rightarrow G^{(0)}$ is a local homeomorphism and that $Z/G$ is compact. Then there exists an integer $n$ and a projection $e\in M_n(C^*(G))$ such that the C*-module $C^*(Z)$ is isomorphic to $e{C^*(G)}^n$.
\end{prop}

\begin{proof} Here, $\alpha_x$ is the counting measure of the discrete fibre $s_Z^{-1}(x)$. We let $r_Z:Z\to Z/G$ be the quotient map.  According to \cite[Theorem 2.8]{mrw:equivalence}, the C*-algebra ${\mathcal K}(C^*(Z))$ of compact operators of the C*-module $C^*(Z)$ is isomorphic to the C*-algebra of the imprimitivity groupoid $(Z*Z)/G$.  Our assumptions imply that this groupoid is \'etale and has a compact unit space. Therefore, its C*-algebra has a unit, namely the characteristic function ${\bf 1}_{Z/G}$ of its unit space.  By density of the subalgebra of finite rank operators, there exist $n\in\N$ and $h_1,\ldots,h_n\in C^*(Z)$ such that ${\rm Id}=\sum_{i=1}^n |h_i\>\<h_i|$. One can take the projection $e=(\<h_i,h_j\>)\in M_n(C^*(G))$.\\
To construct explicitly a projection $e\in M_n(C^*(G))$ giving the module $C^*(Z)$, we find $h_1,\dots,h_n\in C_c(Z)$ satisfying the equality 
\[\sum_{i=1}^n \int_G h_i(z\gamma) \overline{h_i(z'\gamma)} d\lambda^{s_Z(z)}(\gamma)= \begin{cases}
  1 & \text{if } z=z' \\
 0  & \text{otherwise} 
\end{cases}\] 
We proceed as follows. We choose relatively compact open subsets $V_1,\ldots,V_n$ of $Z$ such that for all $i$, the restriction $s_{Z|V_i}$ is one-to-one and  $(r_Z(V_1),\ldots, r_Z(V_n))$ covers $Z/G$. We choose a partition of unity $(b_1,\ldots, b_n)$ subordinate to this open cover. For all $i$, there exists $h_i\in C_c(Z)$ whose support is contained in $V_i$ and such that $\int_G |h_i(z\gamma)|^2d\lambda^{s_Z(z)}(\gamma)=b_i(r_Z(z))$. The projection we are looking for is $e=(<h_i,h_j>)$.
\end{proof}
\begin{rem}\label{n=1}
In particular, if there exists a relatively compact open subset $V$ of $Z$ such that $s_{Z|V}$ is one-to-one and
$r_Z(V)=Z/G$, one obtains a projection $e=\<h,h\>$ in $C^*(G)$ such that $C^*(Z)$ is isomorphic to $e{C^*(G)}$.
\end{rem}
In the next proposition, we use some elements of transverse measure theory (see \cite{con:nci} or \cite[Appendix A.1]{adr:amenable}). Given a locally compact groupoid with Haar system $(G,\lambda)$, a measure $\mu$ on $G^{(0)}$ which is invariant in the sense that $\mu\circ\lambda=(\mu\circ\lambda)^{-1}$ assigns to each transverse function $(Z,\alpha)$ a measure $\nu=\Lambda(\alpha)$ on $Z/G$ in a coherent fashion. This measure is uniquely determined by the relation $\nu\circ\lambda_1=\mu\circ\alpha$, where $\lambda_1$ is the system of measures along the quotient map $Z\to Z/G$ deduced from $\lambda$. Explicitly, $\lambda_1^{[z]}=z\lambda^{s_Z(z)}$ for all $z\in Z$ . On the other hand, the invariant measure $\mu$ defines a trace $\tau$ on $C^*(G)$ such that $\tau(f)=\int_{G^{(0)}}f(x)d\mu(x)$ for all $f\in C_c(G)$. We denote by $\tau_*$ the induced map on $K_0(C^*(G))$ and by $[e]$ the class of the projection $e$ in $K_0(C^*(G))$.

\begin{prop}\cite[Proposition 5.5.6]{ren:rio}\label{volume}. Under the assumptions of the previous proposition, assume that $\mu$ is an invariant measure for $G$. Let $\tau$ be the trace on $C^*(G)$ defined by $\mu$. Then $\tau_*([e])=\nu(Z/G)$, where $\nu$ is the unique measure on $Z/G$ such that $\nu\circ\lambda_1=\mu\circ\alpha$.
\end{prop}

\begin{proof} With $(h_i)$ and $e$ as in the proof of the proposition, we have
\begin{align*}
\tau_*([e])=&\sum_{i=1}^n\tau(\<h_i,h_i\>)=\sum_{i=1}^n \int_{G^{(0)}}\<h_i,h_i\>(x)d\mu(x)\\
=&\sum_{i=1}^n \int_{G^{(0)}}\int_Z |h_i(z)|^2d\alpha_xd\mu(x)=\sum_{i=1}^n\int_Z |h_i|^2d(\mu\circ\alpha)\\
=&\sum_{i=1}^n\int_Z |h_i|^2d(\nu\circ\lambda_1)=\sum_{i=1}^n\int_{Z/G}\int_G |h_i(z\gamma)|^2d\lambda^{s_Z(z)}(\gamma)d\nu(r_Z(z))\\
=&\int_{Z/G}\big(\sum_{i=1}^n\int_G |h_i(z\gamma)|^2d\lambda^{s_Z(z)}(\gamma)\big)d\nu([z])=\int_{Z/G}d\nu([z])=\nu(Z/G)\\
\end{align*}
\end{proof}

\end{document}